\documentclass[reqno,twoside,final]{amsart}
\usepackage{amssymb,amsmath}

\newtheorem{Thm}{Theorem}
\newtheorem{Cor}{Corollary}
\newtheorem{Lem}{Lemma}
\newtheorem{Prop}{Proposition}
\newtheorem{Prob}{Problem}
\theoremstyle{definition}
\newtheorem{Def}{Definition}
\theoremstyle{remark}
\newtheorem{Rem}{Remark}
\newtheorem{Ex}{Example}

\numberwithin{equation}{section}

\newcommand{\thmref}[1]{Theorem~\ref{#1}}
\newcommand{\secref}[1]{\paragraph\ref{#1}}
\newcommand{\lemref}[1]{Lemma~\ref{#1}}
\newcommand{\propref}[1]{Proposition~\ref{#1}}
\newcommand{\probref}[1]{Problem~\ref{#1}}

\let\polishL=\L

\def\rom#1{\/{\rm #1}}

\def\e{\varepsilon}

\def\const{\operatorname{const}}
\def\R{{\mathbb R}}
\def\C{{\mathbb C}}

\let\paragraph=\S
\let\hataccent=\^
\def\S{\varSigma}
\def\Mat{\operatorname{Mat}}

\def\({\left(}
\def\){\right)}
\def\[{\left[}
\def\]{\right]}

\def\H{{\^H}}
\def\diam{\operatorname{diam}}

\def\k{\operatorname{\Bbbk\kern1pt}}

\def\:{\colon}
\def\~#1{\widetilde{#1}}
\def\^#1{\widehat{#1}}
\def\=#1{\check{#1}}
\let\le=\leqslant
\let\ge=\geqslant
\let\ssm=\smallsetminus
\def\<{\left<}
\def\>{\right>}

\begin{document}
\title{Redundant Picard--Fuchs system for Abelian integrals}
\author[D. Novikov, S. Yakovenko]{D. Novikov$^*$, S. Yakovenko$^{**}$}
 \thanks{\leavevmode\hbox to 0pt{\hss$^*$ }TMR fellow in 1998--1999}
 \thanks{\leavevmode\hbox to 0pt{\hss$^{**}$ }Iberdrola professor at
 Universidad de Valladolid in 1999}
\date{January 19, 2000}
\address{Laboratoire de Topologie, Universit\'e de Bourgogne,
Dijon, France}
 \curraddr{Department of Mathematics, Toronto University, Toronto, Canada}
 \email{dmitry@math.toronto.edu}
\address{Department of Mathematics, The Weizmann Institute of
Science, Rehovot, Israel}
 \email{yakov@wisdom.weizmann.ac.il}



\begin{abstract}
We derive an explicit system of Picard--Fuchs differential
equations satisfied by Abelian integrals of monomial forms and
majorize its coefficients. A peculiar feature of this construction
is that the system admitting such explicit majorants, appears only
in dimension approximately two times greater than the standard
Picard--Fuchs system.

The result is used to obtain a partial solution to the tangential
Hilbert 16th problem. We establish upper bounds for the number of
zeros of arbitrary Abelian integrals on a positive distance from
the critical locus. Under the additional assumption that the
critical values of the Hamiltonian are distant from each other
(after a proper normalization), we were able to majorize the
number of all (real and complex) zeros.

In the second part of the paper an equivariant formulation of the
above problem is discussed and relationships between spread of
critical values and non-homogeneity of uni- and bivariate complex
polynomials are studied.
\end{abstract}

\maketitle
\tableofcontents

\section{Tangential Hilbert Sixteenth Problem, complete Abelian integrals
and Picard--Fuchs equations}

The main result of this paper is an explicit derivation of the
Picard--Fuchs system of linear ordinary differential equations for
integrals of polynomial 1-forms over level curves of a polynomial
in two variables, regular at infinity.

The explicit character of the construction makes it possible to derive
upper bounds for the coefficients of this system. In turn, application of
the bounded meandering principle \cite{lleida,meandering} to the system of
differential equations with bounded coefficients allows to produce upper
bounds for the number of complex isolated zeros of these integrals on a
positive distance from the ramification locus.

\subsection{Abelian integrals and tangential Hilbert 16th problem}
If $H(x,y)$ is a polynomial in two real variables, called the
Hamiltonian, and $\omega=P(x,y)\,dx+Q(x,y)\,dy$ a real polynomial
1-form, then the problem on limit cycles appearing in the
perturbation of the Hamiltonian equation,
\begin{equation}\label{tangential-16}
  dH+\e\omega=0,\qquad \e\in(\R,0)
\end{equation}
after linearization in $\e$ (whence the adjective ``tangential'')
reduces to the study of {\em complete Abelian integral\/}
\begin{equation}\label{abelian-integral}
  I(t)=I(t;H,\omega)=\oint_{H=t}\omega,
\end{equation}
where the integration is carried over a continuous family of
(real) ovals lying on the level curves $\{H=t\}$.

\begin{Prob}[Tangential Hilbert 16$^{th}$ problem]
Place an upper bound for the number of real zeros of the Abelian
integral $I(t;H,\omega)$ on the maximal natural domain of
definition of this integral, in terms of $\deg H$ and $\deg
\omega=\max(\deg P,\deg Q)+1$.
\end{Prob}

A more natural version appears after complexification. For an
arbitrary complex polynomial $H(x,y)$ having only isolated
critical points, and an arbitrary complex polynomial 1-form
$\omega$, the integral \eqref{abelian-integral} can be extended as
a multivalued analytic function ramified over a finite set of
points (typically consisting of critical values of $H$). The
problem is to place an upper bound for the number of isolated
complex roots of any branch of this function, in terms of $\deg H$
and $\deg\omega$.

\subsection{Abelian integrals and differential equations}
Despite its apparently algebraic character, the tangential Hilbert
problem still resists all attempts to approach it using methods of
algebraic geometry. Almost all progress towards its solution so
far was based on using methods of analytic theory of differential
equations.

In particular, the (existential) general finiteness theorem by
Khovanski\u{\i}--Var\-chen\-ko \cite{khovanskii,varchenko} claims
that for any finite combination of $d=\deg \omega$ and $n=\deg H$
the number of isolated zeros is indeed uniformly bounded over all
forms and all Hamiltonians of the respective degree. One of the
key ingredients of the proof is the so called {\em Pfaffian
elimination\/}, an analog of the intersection theory for varieties
defined by Pfaffian differential equations \cite{fewnomials}.

Another important achievement, an explicit upper bound for the number of
zeros in the {\em elliptic case\/} when $H(x,y)=y^2+p(x)$, $\deg p=3$ and
forms of arbitrary degree, due to G.~Petrov \cite{petrov}, uses the fact
that the {\em elliptic integrals\/} $I_k(t)=\oint x^{k-1}y\,dx$, $k=1,2$,
in this case satisfy an explicit system of linear first order system of
differential equations with rational coefficients. This method was later
generalized for other classes of Hamiltonians whose level curves are
elliptic (i.e., of genus $1$), see
\cite{horozov-iliev,girard-jebrane,zhang} and references therein. The
ultimate achievement in this direction is a theorem by Petrov and
Khovanskii, placing an {\em asymptotically linear\/} in $\deg\omega$ upper
bound for the number of zeros of arbitrary Abelian integrals, with the
constants being uniform over all Hamiltonians of degree $\le n$
(unpublished). However, one of these constants is purely existential: its
dependence on $n$ is totally unknown.

It is important to remark that all the approaches mentioned above,
require a very basic and easily obtainable information concerning
the differential equations (their mere existence, types of
singularities, polynomial or rational form of coefficients, in
some cases their degree).

\subsection{Meandering of integral trajectories}
A different approach suggested in \cite{era-99} consists in an attempt to
apply a very general principle, according to which integral trajectories
of a polynomial vector field (in $\R^n$ or $\C^n$) have a controllable
meandering (sinuosity), \cite{lleida,meandering}. More precisely, if a
curve of known size is a part of an integral trajectory of a polynomial
vector field whose degree and the magnitude of the coefficients are
explicitly bounded from above, then the number of isolated intersections
between this curve and any affine hyperplane in the ambient space can be
explicitly majorized in terms of these data. The bound appears to be very
excessive: it is polynomial in the size of the curve and the magnitude of
the coefficients, but the exponent as the function of the degree and the
dimension of the ambient space, grows as a tower (iterated exponent) of
height 4.

In order to apply this principle to the tangential Hilbert
problem, we consider the curve parameterized by the {\em monomial
integrals\/},
$$
 t\mapsto(I_1(t),\dots,I_N(t)),\qquad I_i(t)=\oint_{H=t}\omega_i,
$$
where $\omega_i$, $i=1,\dots,N$ are all {\em monomial\/} forms of
degree $\le d$. Isolated zeros of the Abelian integral of an
arbitrary polynomial 1-form $\omega=\sum_i c_i\omega_i$ correspond
to isolated intersections of the above curve with the hyperplane
$\sum c_i I_i=0$. If this monomial curve is integral for a system
of polynomial differential equations with {\em explicitly
bounded\/} coefficients, then the bounded meandering principle
would yield a (partial) answer for the tangential Hilbert 16th
problem.

The system of polynomial (in fact, linear) differential equations
can be written explicitly for the case of {\em hyperelliptic
integrals\/} corresponding to the Hamiltonian $H(x,y)=y^2+p(x)$
with an arbitrary univariate {\em potential\/} $p(x)\in\C[x]$, see
\secref{sec:hyperel} below and references therein. Application of
the bounded meandering principle allowed us to prove in
\cite{era-99} that the number of zeros of hyperelliptic integrals
is majorized by a certain tower function depending only on the
degrees of $n=\deg H=\deg p$ and $d=\deg \omega$. (Actually, it
was done under an additional assumption that all critical values
of $p$ are real, but we believe that this restriction is technical
and can be removed).

\subsection{Picard--Fuchs equations and systems of equations}
In order to generalize the construction from \cite{era-99} for the
case of arbitrary (not necessarily hyperelliptic) Hamiltonians it
is necessary, among other things, to write a system of polynomial
differential equations for Abelian integrals and estimate
explicitly the magnitude of its coefficients.

The mere existence of such a system is well known since times of
Riemann if not Gauss. In today's language, the monodromy group of
any form depends only on the Hamiltonian. Denote by $\mu$ the rank
of the first homology group of a typical affine level curve
$\{H=t\}\subset\C^2$. Then for any collection of 1-forms
$\omega_1,\dots,\omega_\mu$ the {\em period matrix\/} $X(t)$ can
be formed, whose entries are integrals of $\omega_i$ over the
cycles $\delta_1(t),\dots,\delta_\mu(t)$ generating the homology.
If the determinant of this matrix if not identically zero, then
$X(t)$ satisfies a linear ordinary differential equation of the
form
\begin{equation}\label{pf-intro}
 \dot X(t)=A(t)X(t),
 \qquad A(\cdot)\in\Mat_{\mu\times \mu}(\C(t)),
\end{equation}
with a rational matrix function $A(t)$. This system of equations
is known under several names, from {\em Gauss--Manin connection\/}
\cite[especially p.~18]{pham} to {\em Picard--Fuchs system\/} (of
linear ordinary differential equations with rational coefficients,
in full). We shall systematically use the last name.

The rank of the first homology can be easily computed: for a
generic Hamiltonian of degree $n+1$ it is equal to $n^2$. The
degree $\deg A(t)$ can be relatively easily determined if the
degrees of the forms $\omega_i$ are known. However, the choice of
the forms $\omega_i$ may also be a difficult problem for some
Hamiltonians. The matrix $A(t)$ apriori may have poles not only in
the ramification points of the Abelian integrals, which leads to
additional difficulties. But worst of all, this topological
approach gives absolutely no control over the magnitude of the
(matrix) coefficients of the rational (matrix) function $A(t)$.

\subsection{Regularity at infinity and Gavrilov theorems}
Part of these problems problems can be resolved. In particular, if
the Hamiltonian is sufficiently regular at infinity, then all
questions concerning the degrees, can be answered.

\begin{Def}
A polynomial $H(x,y)\in\C[x,y]$ of degree $n+1$ is said to be {\em
regular at infinity\/}, if one of the three equivalent conditions
holds:
\begin{enumerate}
 \item its principal homogeneous part $\H$, a homogeneous polynomial
 of degree $n+1$, is a product of $n+1$ pairwise different linear forms;
 \item $\H$ has an isolated critical point (necessarily of multiplicity
 $\mu=n^2$) at the origin $(x,y)=(0,0)$;
 \item the level curve $\{\H=1\}\subset\C^2$ is nonsingular.
\end{enumerate}
\end{Def}

This condition means that after the natural projective
compactification of the $(x,y)$-plane $\C^2$, all ``interesting''
things still happen only in the finite part of the compactified
plane. In particular, for a polynomial regular at infinity:
\begin{enumerate}
 \item all level curves $\{H=t\}$ intersect the infinite line $\C
 P^1_\infty\subset\C P^2$ transversally,
 \item all critical points $\{(x,y)\:\, dH(x,y)=0\}$ are isolated
 and their number is exactly $\mu=n^2$ if counted with multiplicities,
 \item the rank of the first homology of any regular affine level
 curve $\{H=t\}$ is $\mu=n^2$,
 \item the map $H\:\C^2\to\C^1$ is a topological bundle over the
 set of the regular values of $H$, hence the Abelian integrals can
 be ramified only over the critical values of $H$.
\end{enumerate}

In \cite{gavrilov} L.~Gavrilov proved that for polynomials regular at
infinity, the space of Abelian integrals is finitely generated as a
$\C[t]$-module by $\mu$ {\em basic integrals\/} that can be chosen as
integrals of any $\nu$  forms $\omega_i$ of degree $\le 2n$ whose
differentials form the basis of the quotient space $\Lambda^2/d\H\land
\Lambda^1$, where $\Lambda^k$ is the space of polynomial $k$-forms on
$\C^2$. As a corollary, one can prove that the collection of these basic
integrals satisfies a system of equations \eqref{pf-intro} of size
$\mu\times\mu$ with $\mu=n^2$, and place an upper bound for the degree of
the corresponding matrix function $A(t)$. This system is minimal
(irredundant): generically (for Morse Hamiltonians regular at infinity),
all branches of full analytic continuation of an Abelian integral span
exactly $\mu$-dimensional linear space.

From this theorem one can also derive further information
concerning the Picard--Fuchs system. Namely, one can prove that if
in addition to being regular at infinity, $H$ is a Morse function
on $\C^2$, then the matrix $A(t)$ of the Picard--Fuchs system
\eqref{pf-intro} has only simple poles (Fuchsian singularities) at
the critical values of the Hamiltonian and only at them (the point
$t=\infty$ is a regular though in general non-Fuchsian
singularity).

However, these results do not yet allow an explicit majoration of
the coefficients (e.g., the residue matrices) of the matrix
function $A(t)$ in \eqref{pf-intro}.

\subsection{Redundant Picard--Fuchs system: the first main result}
We suggest in this paper a procedure of explicit derivation of the
Picard--Fuchs system of equations, based on the division by the
gradient ideal $\<H_x,H_y\>\subset\C[x,y]$ in the polynomial ring.
It turns out that if instead of choosing $\mu=n^2$ forms of degree
$\le 2n$ constituting a basis modulo the gradient ideal, one takes
all $\nu=n(2n-1)$ cohomologically independent monomial forms of
degree $\le 2n$, then the resulting Picard--Fuchs system can be
written in the form
\begin{equation}\label{tAB-intro}
 (tE-A)\dot X(t)=B X(t),\qquad A,B\in\Mat_{\nu\times\nu}(\C),
\end{equation}
where $E$ is the identity matrix, and $X(t)$ is the rectangular
period $\nu\times \mu$-matrix.

The procedure of deriving the system \eqref{tAB-intro}, being completely
elementary, can be easily analyzed and upper bounds for the matrix norms
$\|A\|$ and $\|B\|$ derived. These bounds depend on the magnitude of the
all non-principal terms $H-\H$ of the Hamiltonian, relative to the
principal part $\H$.

More precisely, we introduce a normalizing condition ({\em
quasimonicity\/}) on the homogeneous part: this condition plays
the same role as the assumption that the leading term has
coefficient $1$ for univariate polynomials. The quasimonicity
condition can be always achieved by an affine change of variables,
provided that $H$ is regular at infinity, hence it is not
restrictive. \thmref{main}, our first main result, allows to place
an upper bound for the norms $\|A\|+\|B\|$ in terms of the norm
(sum of absolute values of all coefficients) of the non-leading
part $H-\H$, assuming that $H$ is quasimonic.

\subsection{Corollaries: theorems on zeros}\label{theoremsonzeros}
The above information on coefficients of Picard--Fuchs system
already suffices to apply the bounded meandering principle and
obtain {\em an explicit upper bound\/} for the number of zeros of
complete Abelian integrals {\em away from the critical locus of
the Hamiltonian\/} (\thmref{onzeros}), which seems to be the first
known explicit result of that kind.

In addition to this bound valid for {\em some\/} zeros and {\em almost
all\/} Hamiltonians, one can apply results (or rather methods) from
\cite{roitman}. If in addition to the quasimonicity and bounded lower
terms, all {\em critical values\/} $t_1,\dots,t_\mu$ of the Hamiltonian
$H$ are far away from each other (i.e., a {\em lower\/} bound for
$|t_i-t_j|$ is known for $i\ne j$), then one can majorize the number of
zeros on any branch of the Abelian integral by a function depending only
on $n,d$ and the minimal distance between critical values. The accurate
formulation is given in \thmref{noncollision}.

\subsection{Equivariant formulation}
However, the description given by \thmref{main}, is not completely
sufficient for further advance towards solution of the tangential Hilbert
problem by studying zeros of Abelian integrals {\em near the critical
locus\/} when the latter (or some part of it) shrinks to one point of high
multiplicity.

One reason is that in order to run an inductive scheme similar to
that constructed in \cite{era-99}, one has to make sure that the
Hamiltonian $H\:\C^2\to\C^1$ can be rescaled (using affine
transformations in the preimage $\C^2$ and the image $\C^1$) so
that {\em simultaneously\/}:
\begin{enumerate}
 \item the critical values of $H$ do not tend to each other (e.g., their
 diameter is bounded from below by $1$), and
 \item the ``non-homogeneous part'' $H-\H$ is bounded by a
 constant explicitly depending on $n$
\end{enumerate}
(each of the two conditions can be obviously satisfied
separately).

Another, intrinsic reason is the equivariance (or, rather
precisely, non-invariance of neither \thmref{main} nor
Theorems~\ref{onzeros} and \ref{noncollision}) by the above affine
group action. In order to be geometrically sound, all assertions
should be related to a certain privileged {\em affine\/} chart on
the $t$-plane. Since our future goal is to study a neighborhood of
the critical locus, it is natural to choose the privileged chart
so that the critical locus will not shrink into one point.

More detailed explanations and motivations are given in
\secref{motivation} below, where we formulate several problems all
in the following sense: for a polynomial whose principal
homogeneous part is normalized (in a certain sense) and whose
critical values are explicitly bounded, it is required to place an
upper bound for the ``non-homogeneous'' part, eventually after a
suitable translation (which does not affect the principal part,
naturally).

\subsection{Geometry of critical values of polynomials}
The reason why several problems of the above type were formulated
instead of just one, is very simple: {\em we do not know a
complete solution\/}, so partial, existential or limit cases were
considered as intermediate steps towards the ultimate goal. In
\secref{spread} we prove that:
\begin{itemize}
 \item if a monic complex polynomial $p(x)=x^{n+1}+\cdots\in\C[x]$ has all
 critical values in the unit disk, then its roots form a point set
 of diameter $<11$ (\thmref{roots}) and
 hence by a suitable translation the
 norm of the non-principal part can be made $\le 12^{n+1}$ (this
 gives a complete solution in the univariate and hyperelliptic cases);
 \item all critical values of a Hamiltonian regular at infinity,
 cannot simultaneously coincide unless the Hamiltonian is essentially
 homogeneous (\thmref{one-val});
 \item for any normalized principal part $\H$ there exists an
 upper bound for $H-\H$ (eventually after a suitable translation),
 provided that the critical values of $H$ are all in the unit
 disk (Corollary to \thmref{one-val}).
\end{itemize}
All these are positive results towards solution of the problem on
critical values. It still remains to compute the upper bound from
the last assertion explicitly: the proof below does not provide
sufficient information for that.

However, it can be shown already in simple examples that this
bound {\em cannot be uniform\/} over all homogeneous parts. As
some of the linear factors from $\H$ approach too closely to each
other, an explosion occurs and the non-principal part may be
arbitrarily large without affecting the ``moderate'' critical
values. The phenomenon can be seen as ``almost occurrence'' of
atypical values, ramification points for Abelian integrals that
are not critical values of $H$: such points are known to appear if
the principal part $\H$ has a non-isolated singularity.

\subsection*{Acknowledgements}
We are grateful to J.-P.~Fran\-\c coise, L.~Gavrilov, Yu.~Ilya\-shen\-ko,
A.~Khovanskii, P.~Milman, R.~Moussu, R.~Roussarie and Y.~Yomdin for
numerous discussions and many useful remarks. Bernard Teissier suggested
an idea that finally developed into the proof of \lemref{homogen} below.
Lucy Moser provided us with a counterexample (\secref{sec:demo1val}
below).

We are grateful to all our colleagues from Laboratoire de Topologie,
Universit\'e de Bourgogne (Dijon) and Departamento de Algebra, Geometr\'\i
a y Topolog\'\i a, Universidad de Valladolid, where a large part of this
work has been done. They made our stays and visits very stimulating.

\section{Picard--Fuchs system in the hyperelliptic case}
\label{sec:hyperel}

\subsection{Gelfand--Leray residue}
The derivative of an Abelian integral $\oint_{H=t}\omega$ can be
computed as the integral over the same curve of another 1-form
$\theta$ called the Gelfand--Leray derivative (residue). More
precisely, if a pair of polynomial  1-forms $\omega,\theta$
satisfies the identity $d\omega=dH\land \theta$, then for any
continuous family of cycles $\delta(t)$ on the level curves
$\{H=t\}$
\begin{equation}\label{gelfand-leray}
  \frac d{dt}\oint_{\delta(t)}\omega=\oint_{\delta(t)}\theta,
  \qquad\forall \delta(t)\subset\{H=t\}
\end{equation}
(the Gelfand--Leray formula). The identity remains true if
$\theta$ is only meromorphic but has zero residues after
restriction on each curve $H=t$.

The identity between $\omega$, $dH$ and $\theta$ explains the
standard notation $\theta=d\omega/dH$: to find $\theta$, one has
to divide $d\omega$ by $dH$. In general this division is not
possible in the class of polynomial 1-forms, but one can always
divide $d\omega$ by $dH$ {\em with remainder\/}: the corresponding
identity after integration will give a differential equation
relating Abelian integrals with their derivatives.

We illustrate this idea by deriving explicitly the Picard--Fuchs
system for hyperelliptic Hamiltonians. In the hyperelliptic case
the outlined approach yields a complete and in some sense minimal
(irredundant) system that could be in principle derived by a
number of different ways, e.g., as in \cite{givental}. Moreover,
using the explicit nature of Euclid's algorithm of division of
univariate polynomials, one can produce explicit upper bounds for
the magnitude of the coefficients of the resulting equations, that
are difficult (if possible at all) to obtain applying methods from
\cite{givental}. The constructions from this section serve as a
paradigm for further exposition \secref{redundant}.

\subsection{Division by polynomial ideals and 1-forms}
Let $q_1,q_2\in\C[x,y]$ be a pair of polynomials generating the
ideal $\<q_1,q_2\>\subset\C[x,y]$ that has a finite codimension
$\mu$. By definition, this means that there exist $\mu$
polynomials $r_1,\dots,r_\mu\in\C[x,y]$ (the remainders) such that
any polynomial $f\in\C[x,y]$ admits representation $v=q_1 u_2-q_2
u_1+\sum_1^\mu \lambda_i r_i$ with  polynomials
$u_1,u_2\in\C[x,y]$ and constants $\lambda_i\in\C$.

It is convenient to interpret this identity as a division formula
for polynomial 2-forms: any polynomial 2-form
$\Omega=f(x,y)\,dx\land dy$ can be divided by the given 1-form
$\xi=q_1\,dx+q_2\,dy$ with the ``incomplete ratio''
$\eta=u_1\,dx+u_2\,dy$ and the remainder that is a linear
combination of the 2-forms $\Omega_i=r_i\,dx\land dy$,
$$
 \Omega=\xi\land \eta+\sum_{i=1}^\mu \lambda_i \Omega_i.
$$
Denoting by $\Lambda^k$, $k=0,1,2$, the modules (over the ring
$\C[x,y]$) of polynomial $k$-forms on $\C^2$, we say that the
tuple of 2-forms $\{\Omega_i\}_1^\mu$ generates the quotient
$\Lambda^2/\xi\land \Lambda^1$.

The {\em gradient ideal\/} $\<H_x,H_y\>$  has a finite codimension
provided that $H$ has only isolated critical points. In this case
we will usually apply the division formula to a differential
$\Omega=d\omega$ of a polynomial 1-form and write the generators
explicitly as $\Omega_i=d\omega_i$ for appropriate polynomial
primitives $\omega_i\in\Lambda^1$:
\begin{equation}\label{div2by1}
  d\omega=dH\land \eta+\sum_{i=1}^\mu\lambda_i\,d\omega_i.
\end{equation}
This means that the Gelfand--Leray derivative of the form
$\omega-\sum_1^\mu\lambda_i\omega_i$ can be found in the class of
{\em polynomial\/} 1-forms, $\eta\in\Lambda^1$.

\subsection{Derivation of the Picard--Fuchs system in the hyperelliptic
case} Throughout this section we assume that
$H(x,y)=\frac12y^2+p(x)$, where $p\in\C[x]$ is a {\em monic\/}
polynomial of degree $n+1$ in one variable:
$p(x)=x^{n+1}+\sum_{i=0}^{n-1}c_ix^i$, $\sum|c_i|=c$.

The gradient ideal and the corresponding quotient in this case can
be easily computed:
\begin{equation*}
  \<H_x,H_y\>=\<p'(x),y\>, \quad
  \C[x,y]/\<H_x,H_y\>\simeq\C[x]/\<x^n\>
  \simeq\bigoplus_{k=1}^n \C x^{k-1},
\end{equation*}
so that the quotient algebra is an algebra of truncated univariate
polynomials of degree $\le n-1$. This observation motivates the
following computation.

Denote by $\omega_i=x^{i-1}y\,dx$, $i=1,\dots,n$, the differential
1-forms whose derivatives $d\omega_i=x^{i-1}\,dx\land dy$ generate
$\Lambda^2/dH\land \Lambda^1$. Then
\begin{equation*}
\begin{aligned}
 H\,d\omega_i&=\(\tfrac12y^2+p(x)\)x^{i-1}\,dx\land dy
 \\
 &=\[\tfrac12x^{i-1}y H_y+\(b_i(x)H_x+a_i(x)\)\]\,dx\land dy
 \\
 &=\(\tfrac12x^{i-1}y\,dx-b_i(x)\,dy\)\land dH+a_i(x)\,dx\land dy
 \\
 &=\biggl[\tfrac12\omega_i+\sum_{j=1}^n
 b_{ij}\omega_j+d(yb_j(x))\biggr]\land dH+\sum_{j=1}^n a_{ij}d\omega_j,
\end{aligned}
\end{equation*}
where we used the following identities:
\begin{enumerate}
\item[(i)] division with remainder: the polynomial $x^{i-1}p(x)$ of
degree $n+i$ is divided out by $p'(x)=H_x$ as
\begin{equation}\label{div-rem-1}
  x^{i-1}p(x)=b_i(x)p'(x)+a_i(x),\qquad \deg b_i\le i,\ \deg a_i\le n,
\end{equation}
\item[(ii)] the form $b_i(x)\,dy$ is represented as a linear combination
of the basic forms modulo an exact term:
\begin{equation}\label{modulo-exact}
  b_i(x)\,dy=d(yb_i)-b_i'(x)y\,dx=\sum_{j=1}^{i} b_{ij}x^{j-1}y\,dx+dF_i,
\end{equation}
since the degree of $b_i'\in\C[x]$ never exceeds $i-1$;
\item[(iii)] the remainders $a_i(x)\,dx\land dy$ can be represented as
linear combinations of $d\omega_j$:
\begin{equation}\label{remainder}
  a_i(x)\, dx\land dy=\sum_{j=1}^n a_{ij}x^{j-1}\,dx\land
  dy=\sum_{j=1}^n a_{ij}\,d\omega_j.
\end{equation}
\end{enumerate}
Integrating over closed ovals of the level curves $H=t$ (so that
the exact forms $dF_i$ disappear) and using the Gelfand--Leray
formula \eqref{gelfand-leray}, we conclude with the system of
linear ordinary differential equations
\begin{equation}\label{pf-hyperell}
  t\dot I_i-\sum_{j=1}^n a_{ij}\dot I_j=\tfrac12 I_i+\sum_{j=1}^n
  b_{ij}I_j
\end{equation}
or, in the matrix form,
\begin{equation}\label{pf-h-matrix}
  (t-A)\dot I=BI,\qquad I\in\C^n,\ A,B\in\Mat_{n\times n}(\C),
\end{equation}
where, obviously, $I_j(t)=\oint\omega_j$ are the Abelian integrals
and $I=(I_1,\dots,I_n)$ the column vector.

\begin{Rem}
The computation above does not depend on the choice of the cycle of
integration, therefore the system of equations will remain valid if we
replace the column vector $I$ by the {\em period matrix\/} $X(t)$ obtained
by integrating all forms $\omega_i$ over all {\em vanishing cycles\/}
$\delta_j(t)$, $j=1,\dots,n$ (see \cite{avg}) on the hyperelliptic level
curves.
\end{Rem}

\begin{small}
The matrices $A,B$ can be completely described using the division
process \eqref{div-rem-1}. Let $x_*\in\C$ be a critical point of
$p$ and $t_*=p(x_*)$ the corresponding critical value. Then
\eqref{div-rem-1} imply that the column vector
$(1,x_*,x_*^2,\dots,x_*^{n-1})\in\C^n$ is the eigenvector of $A$
with the eigenvalue $t_*$, which gives a complete description
(eigenbasis and eigenvalues) of $A$.

Entries of the matrix $B$ can be described similarly: $b_{ij}=0$
for $j>i$ because of the assertion about degrees of $b_j(x)$, so
$B$ is triangular. The diagonal entries can be easily computed by
looking at the leading terms: since $p$ is monic,
$b_i(x)=\frac{x^i}{n+1}+\cdots$, hence
$-b_i'(x)=-\frac{i}{n+1}x^{i-1}+\cdots$. Finally, the diagonal
entries $\{-\frac i{n+1}+\frac12\}_{i=1}^n$ form the spectrum of
$B$.

However, knowledge of the critical values of $H$ is not yet
sufficient to produce an upper bound for the norms $\|A\|,\|B\|$,
since the conjugacy by the Vandermonde matrix (whose columns are
the above eigenvectors $(1,x_j,x_j^2,\dots,x_j^{n-1})^T$,
$j=1,\dots,n$) may increase arbitrarily the norm of the diagonal
matrix $\operatorname{diag}(t_1,\dots,t_n)$, where $t_j$ are all
critical values of $H$ (or $p$, what is the same). On the
contrary, a linear change in the space of 1-forms that makes $A$
diagonal, can increase in an uncontrollable way the norm of the
matrix $B$, whose eigenbasis differs from the standard one by a
triangular transformation. It is the explicit division procedure
that allows to majorize the matrix norms.\par
\end{small}

\subsection{Bounds for the matrix norms}\label{bounds-1-dim}
For a polynomial $p\in\C[x]$ let $\|p\|$ be the sum of absolute
values of its coefficients (sometimes it is called the length of
$p$). It has the advantage of being multiplicative,
$\|pq\|\le\|p\|\cdot\|q\|$.

\begin{Prop}\label{bounds-univariate-division}
If $q=x^n+\cdots\in\C[x]$ is a monic polynomial with
$\|q-x^n\|=c$, then any other polynomial $f\in\C[x]$ of degree
$d\ge n$ can be divided with remainder,
\begin{equation}\label{div-norm-1}
  f(x)=b(x)q(x)+a(x),\qquad \deg a\le n-1,
\end{equation}
so that
\begin{equation}\label{div-norm-2}
  \|b\|+\|a\|\le K\|f\|,\qquad K=1+C+C^2+\cdots+C^{d-n},\quad
  C=1+c=\|q\|.
\end{equation}
\end{Prop}

\begin{proof}
The proof goes by direct inspection of the Euclid algorithm of
univariate polynomial division. The assertion of the Proposition
is trivial for $q=x^n$: in this case the string of coefficients of
$r$ has to be split into two, and immediately we have the
decomposition $r=b x^n+a$ with $\|b\|+\|a\|=\|r\|$.

The general nonhomogeneous case is treated by induction. Suppose
that the inequality \eqref{div-norm-2} is valid for any polynomial
$\~f$ of degree $\le d-1$ (for $d=n-1$ it is trivially satisfied
by letting $b=0$ and $a=\~f$). Take a polynomial $f$ of degree $d$
and write the identity $f=bx^n+a=b q+b(x^n-q)+a=bq+\~f$, where the
polynomial $\~f=a+b(x^n-q)$ is of degree $\le d-1$ and has the
norm explicitly bounded: $\|\~f\|\le
c\|b\|+\|a\|\le(1+c)(\|b\|+\|a\|)\le C\|f\|$. By the induction
assumption, $\~f$ can be divided, $\~f=\~bq+\~a$, with the norms
satisfying the inequality \eqref{div-norm-2}. Collecting
everything together, we have $f=(b+\~b)q+\~a$ and
$\|b+\~b\|+\|\~a\|\le
\|f\|+C\|f\|(1+C+\cdots+C^{d-1-n})\le\|f\|(1+\cdots+C^{d-n})$.
\end{proof}

As a corollary to this Proposition and the explicit procedure of the
division, we obtain upper bounds for norms of the matrices $A,B$. Recall
that we use the $\ell^1$-norm on the ``space of columns'', so the norm of
a matrix $A=\(a_{ij}\)_{i,j=1}^n$ is
\begin{equation}\label{matrix-norm}
  \|A\|=\max_{j=1,\dots,n}\sum_{i=1}^n|a_{ij}|
\end{equation}

\begin{Thm}\label{hyperelliptic-AB-bounds}
Suppose that $p(x)=x^{n+1}+\sum_{i=0}^{n-1}c_i x^i$ is a monic
polynomial of degree $n+1$ and the non-principal part of $p$ is
explicitly bounded: $\sum_{i=0}^{n-1}|c_i|\le c$.

Then the entries of the matrices $A,B$ determining the
Picard--Fuchs system \eqref{pf-h-matrix} are explicitly bounded:
\begin{equation}\label{normsAB}
  \|A\|+\|B\|\le n^2(1+C+\cdots+C^{n+1}),\qquad C=1+c=\|p\|.
\end{equation}
\end{Thm}

\begin{proof}
The derivative $p'(x)$ is not monic, but the leading coefficient
is explicitly known: $p'(x)=(n+1)(x^n+\cdots)\|$, with the
non-principal part denoted by the dots bounded by $c$ in the sense
of the norm. Applying \propref{bounds-univariate-division} to
$q=p'/(n+1)$, we see that any polynomial can be divided by $p'$
and the same inequalities \eqref{div-norm-2} would hold (since
$n+1\ge 1$).

Thus we have $\|b_i\|+\|a_i\|\le K\|x^{i-1}\|\|p\|=KC$, where
$K=1+C+\cdots+C^n$, then obviously $\|b_j'\|\le n\|b_j\|$ and
finally for the sum of matrix elements $A,B$ occurring in the
$i$th line, we produce an upper bound $\tfrac12+\sum_j
|b_{ij}|+\sum_j |a_{ij}|\le n C(1+C+\cdots+C^n)+\tfrac12\le
n(1+C+\cdots+C^{n+1})$. Clearly, this means that every entry of
these matrices is majorized by the same expression and therefore
for the matrix $\ell^1$-norms on $\C^n$ we have the required
estimate.
\end{proof}

\subsection{Digression: doubly hyperelliptic Hamiltonians}
The algorithm suggest\-ed above, works with only minor
modifications for {\em doubly hyperelliptic\/} Hamiltonians having
the form $H(x,y)=p(x)+q(y)$ (the hyperelliptic case corresponds to
$q(y)=\frac12 y^2$). Assume that $n+1=\deg_x p$, $m+1=\deg_y q$
(there is no reason to require that $n=m$).

In this case the quotient algebra by the gradient ideal is
generated by $nm$ monomials $x^i y^j$, $0\le i\le n-1$, $0\le j\le
m-1$. We claim that any collection of monomial primitives
$\omega_{ij}$ to the monomial 2-forms
$d\omega_{ij}=x^iy^j\,dx\land dy$ satisfies a system of $nm$
equations having the same form \eqref{pf-h-matrix} though a
different size. Indeed,
\begin{equation*}
  H\,d\omega_{ij}=p(x)x^i\,dx\land y^j\,dy-q(y)y^j\,dy\land
  x^i\,dx.
\end{equation*}
Dividing the 1-form $p(x)x^i\,dx$ with remainder by the 1-form
$dp(x)$, we express the former as $b_i(x)\,dp(x)+a_i(x)\,dx$ with
$\deg b_i\le i+1\le n+1$, $\deg a_i\le n-1$ and multiply the
result by $y^j\,dy$. The second term can similarly be rewritten
involving the representation $q(y)y^j=b^*_j(y)\,dq+a_j^*(y)\,dy$.
Putting everything together, we conclude that
$$
 H\,d\omega_{ij}
 =[dp(x)\land b_i(x)y^j\,dy
 -dq(y)\land b_j^*(y)x^i\,dx]+[a_i(x)y^j-x^ia_j^*(y)]\,dx\land dy
$$
Since $dH=dp(x)+dq(y)$, we see that the first bracket is actually
the wedge product $dH\land \eta_{ij}$, where
$\eta_{ij}=b_i(x)y^j\,dy+b_j^*(y)x^i\,dx$ is a polynomial 1-form
whose differential
$$
 d\eta_{ij}=\biggl(\frac{\partial b_i}{\partial x}y^j-
 \frac{\partial b_j^*}{\partial y}x^i\biggr)\,dx\land dy
$$
has the coefficient of degree $\le i$ in $x$ and $\le j$ in $y$
and hence can be expanded as a linear combination of the forms
$\omega_{ij}$ modulo an exact form. The second bracket, being a
2-form with coefficient of degrees $\le n-1$ in $x$ and $\le m-1$
in $y$, is a linear combination of the forms $d\omega_{ij}$. Thus
we have the equations
$$
\begin{gathered}
 H\,d\omega_{ij}=dH\land
 \biggl(\sum_{k,l=0}B_{ij,kl}\,\omega_{kl}+dF_{ij}\biggr)+
 \sum_{k,l}A_{ij,kl}\,d\omega_{kl},
 \\
 i,k=0,\dots,n-1,\quad j,l=0,\dots,m-1.
 \end{gathered}
$$

Rescaling $x$ and $y$ by appropriate factors {\em
independently\/}, we can assume that the polynomials $p(x)$ and
$q(y)$ are both monic. Then all divisions will be bounded provided
that the norms $\|p\|$ and $\|q\|$ are explicitly bounded, and in
a way completely similar to the arguments from
\secref{bounds-1-dim}, we can derive upper bounds for the matrix
coefficients $A_{ij,kl},B_{ij,kl}$.

Thus the case of doubly hyperelliptic Hamiltonians does not differ
much from the ordinary hyperelliptic case, at least as far as the
Picard--Fuchs systems for Abelian integrals are concerned.

\subsection{Discussion}
The Picard--Fuchs system written in the form \eqref{pf-h-matrix}
for a generic  hyperelliptic Hamiltonian (with the potential
$p(x)$ being a Morse function on $\C$), is a system remarkable for
several instances:
\begin{itemize}
 \item it possesses only Fuchsian singularities (simple poles)
 both at all finite singularities $t=t_j$, $j=1,\dots,n$, and at
 infinity;
 \item it has no apparent singularities: all points $t_j$ are
 ramification points for the fundamental system of solutions $X(t)$
 that is obtained by integrating all forms $\omega_i$ over all
 vanishing cycles $\delta_j(t)$ (the period matrix);
 \item it is minimal in the sense that
 analytic continuations of any column of the period matrix $X(t)$
 along all closed loops span the entire space $\C^n$.
 \item its coefficients can be explicitly bounded in terms of
 $\|H\|$.
\end{itemize}
(All these observations equally apply to doubly hyperelliptic
Hamiltonians.)

In the next section we generalize this result for arbitrary
bivariate Hamiltonians. It will be impossible to preserve all
properties, and we shall concentrate on the derivation of {\em
redundant\/} system, eventually exhibiting apparent singularities,
but all of them (including that at infinity) Fuchsian and with
explicitly bounded coefficients.

\section{Derivation of the redundant Picard--Fuchs
system}\label{redundant}

\subsection{Notations and conventions}
Recall that $\Lambda^k$ denote the spaces of polynomial $k$-forms
on $\C^2$ for $k=0,1,2$. They will be always equipped with the
$\ell^1$-norms: the norm of a form is always equal to the sum of
absolute values of all its coefficients. This norm behaves
naturally with respect to the (wedge) product: for any two forms
$\eta\in\Lambda^k$, $\theta\in\Lambda^l$, $0\le k+l\le 2$, we
always have $\|\eta\land\theta\|\le\|\eta\|\cdot\|\theta\|$.

It is also convenient to grade the spaces of polynomial forms so
that the degree of a $k$-form is the maximal degree of its
(polynomial) coefficients plus $k$. Under this convention the
exterior derivation is degree-preserving: $\deg d\theta=\deg
\theta$ (unless $d\theta=0$). An easy computation shows that
$\|d\theta\|\le\deg\theta\cdot\|\theta\|$ for any 0- and 1-form
$\theta$. On several occasions the finite-dimensional linear space
of $k$-forms of degree $\le d$ will be denoted by $\Lambda^k_d$.

If $\omega\in\Lambda^1$ is a polynomial 1-form and
$H\in\Lambda^0$, then by $d\omega/dH$ is always denoted the
Gelfand--Leray derivative \eqref{gelfand-leray}, while by
$\frac{d\omega}{dx\land dy}$ we denote the polynomial coefficient
of the 2-form $d\omega$.

The space $\Lambda^2$ sometimes will be identified with
$\Lambda^0\simeq\C[x,y]$, the submodule $dH\land\Lambda^1$ with
the gradient ideal $\<H_x,H_y\>\subset\C[x,y]$, and the local
algebra as a linear space over $\C$ with the quotient
$\Lambda^2/dH\land \Lambda^1$.

\subsection{Normalizing conditions and quasimonic Hamiltonians}
In the ring $\C[x]$ of univariate polynomials division by the
principal ideal $\<p\>$ is a linear operator whose norm can be
controlled in terms of $\|p\|$ provided that the leading term of
$p$ is bounded from below, in particular when the polynomial is
monic. The definition below introduces a generalization of this
condition for ideals in the ring $\C[x,y]$ of bivariate
polynomials.

Recall that two {\em homogeneous\/} polynomials $a,b\in\C[x,y]$ of
the same degree $n$ have no common linear factors if and only if
their resultant is nonzero and hence the Sylvester matrix is
invertible. In this case an arbitrary homogeneous polynomial $f$
of degree $2n-1$ can be represented as $f=au+bv$ with appropriate
(uniquely defined) homogeneous polynomials $u,v$ of degree $n-1$
each.

\begin{Def}
A pair of homogeneous polynomials $a,b\in \C[x,y]$ of degree $n$
is said to be {\em normalized\/} if the linear operator
$(u,v)\mapsto au+bv$ restricted on the subspace of pairs of
homogeneous polynomials of degree $n-1$, has the inverse of the
unit norm, in other words, if any homogeneous polynomial $f$ of
degree $2n-1$ can be represented as $f=au+bv$ with an explicit
control over norms of  the homogeneous ``ratios'' $u,v$ of degree
$n-1$:
\begin{equation}\label{bdp}
  f=au+bv,\qquad \|u\|+\|v\|\le\|f\|.
\end{equation}
\end{Def}

\begin{Def}
A homogeneous polynomial 1-form $\eta=a\,dx+b\,dy$ of degree $n+1$
is normalized if its coefficients $a,b\in\C[x,y]$ form a
normalized pair.
\end{Def}

For nonhomogeneous objects we impose normalizing conditions on
their principal homogeneous part.

\begin{Def}
A polynomial 1-form $\xi\in \Lambda^1$ of degree $n$ is {\em
normalized at infinity\/}, if its principal homogeneous part
$\^\xi$ is normalized.

A Hamiltonian $H(x,y)\in\C[x,y]$ of degree $n+1$ is said to be
normalized at infinity or {\em quasimonic\/}, if $dH$ is
normalized at infinity in the sense of the previous definition.
\end{Def}

\begin{Rem}
To be normalized at infinity has nothing to do with the
$\ell^1$-norm of a form or Hamiltonian. We will mostly use the
term ``quasimonic''.
\end{Rem}

\subsection{Balanced Hamiltonians}
In order to simplify the calculations below, we impose additional
normalizing condition on $H$ meaning that the non-principal (low
degree) terms are not dominating the principal part.

\begin{Def}
A Hamiltonian $H\in\C[x,y]$ will be called {\em balanced\/}, if it
is quasimonic (the principal homogeneous part $\^H$ is normalized)
and $\|H-\H\|\le 1$.
\end{Def}

For a balanced Hamiltonian, its differential $dH$ is a 1-form that
is (by definition) normalized at infinity and differs from its
principal homogeneous part $d\H$ by the form of degree $n$ and
$\|dH-d\H\|\le n$.

The two conditions, normalization at infinity and that of balance
between principal and non-principal parts, can be obtained
simultaneously by suitable affine transformations. If the
Hamiltonian $H$ is regular at infinity, then after a suitable
choice of $\lambda\in\C$ one can make any of the two polynomials,
$\lambda H(x,y)$ or $H(\lambda x,\lambda y)$ being normalized at
infinity (the same refers to 1-forms). Furthermore, if $H$ is
already quasimonic, one can always choose a suitable
$\lambda\in\C$ so that
$\lambda^{n+1}H(\lambda^{-1}x,\lambda^{-1}y)$ will be balanced
while remaining quasimonic.

\subsection{Lemma on bounded division}
Division by a balanced 1-form is a linear operator whose norm can
be easily controlled.

Let $\xi\in\Lambda^1$ be a polynomial 1-form of degree $n+1$
normalized at infinity, with the principal homogeneous part
denoted by $\^\xi$.

\begin{Lem}\label{bd}
Any polynomial 2-form $\Omega\in\Lambda^2$ can be divided with
remainder by $\xi$,
\begin{equation}\label{div-rem-form}
  \Omega=\xi\land\eta+\Theta,
\end{equation}
where the remainder $\Theta\in\Lambda^2$ is a 2-form of degree
$\le 2n$ and the ``incomplete ratio'' $\eta\in\Lambda^1$ is a
1-form of degree $\deg\Omega-\deg\xi$.

The decomposition \eqref{div-rem-form} is in general non-unique.
However, one can always find $\eta$ and $\Theta$ so that if
$\|\xi-\^\xi\|=c$, then
\begin{equation}\label{div-rem-bounds}
  \|\eta\|+\|\Theta\|\le K\|\Omega\|,
  \qquad K=(1+C+\cdots+C^{d-2n}),
\end{equation}
where $C=c+1$ and $d=\deg\Omega$.
\end{Lem}

\begin{proof}
The proof reproduces almost literally the division algorithm for
univariate polynomials, see
Proposition~\ref{bounds-univariate-division}.

1. For a homogeneous form $\Omega=f\,dx\land dy$ of degree $2n+1$
the divisibility $\Omega=\^\xi\land\eta$ by the homogeneous form
$\^\xi=a\,dx+b\,dy$ is the same as the representation \eqref{bdp}
(recall that our convention concerning the degrees of the form
means that in this case $\deg f=2n-1$). From the normalization
condition it follows then that
$\|\eta\|=\|u\|+\|v\|\le\|f\|=\|\Omega\|$ simply by definition.

2. Writing the division identities for all monomial forms of
degree $2n+1$, multiplying them by arbitrary monomials and adding
results we see then that any polynomial $2$-form $\Omega$ {\em
containing no terms of degree $2n$ and less\/}, can be divided by
$\^\xi$ and the norm of the ``ratio'' $\~\eta$ does not exceed
$\|\Omega\|$. Finally, {\em any\/} form can be represented as the
sum of a ``remainder'' $\~\Theta$, the collection of terms of
degree $\le 2n$, and the higher terms divisible by $\^\xi$.

All together this means that if $\^\xi$ is a homogeneous
normalized 1-form of degree $n+1$, then any polynomial 2-form
$\Omega$ can be divided out as
\begin{equation}\label{div-by-homogen}
  \Omega=\^\xi\land\~\eta+\~\Theta,\qquad
  \|\~\eta\|+\|\~\Theta\|\le\|\Omega\|,
  \qquad\deg\~\eta\le\deg\Omega-\deg\^\xi.
\end{equation}

3. To divide by a nonhomogeneous form $\xi$ normalized at
infinity, we first divide by its principal part $\^\xi$ as in
\eqref{div-by-homogen}. Then
\begin{equation}\label{div-rem-ind}
  \Omega=\xi\land\~\eta+(\^\xi-\xi)\land\~\eta+\~\Theta=
  \xi\land\~\eta+\~\Omega.
\end{equation}
It remains to notice that
$\|\~\eta\|\le\|\~\eta\|+\|\~\theta\|\le\|\Omega\|$ and $\~\Omega$
is a new 2-form whose degree is strictly less than $d=\deg
\Omega$, provided that $d>2n$. Since the norm of $\xi-\^\xi$ is
explicitly bounded by $c$, we have
$$
\|\~\Omega\|\le c\|\~\eta\|+\|\~\Theta\|\le
(1+c)(\|\~\eta\|+\|\~\Theta\|)\le C\|\Omega\|.
$$
We may now continue by induction, accumulating the divided parts
$\~\eta$ and reducing the degrees of ``incomplete remainders''
$\~\Omega$ until the latter become less or equal to $2n$. More
accurately, we use the inductive assumption to divide out
$\~\Omega=\xi\land\eta'+\Theta$ with
$\|\eta'\|+\|\Theta\|\le\|\~\Omega\|(1+C+\cdots+C^{d-1-2n})\le
\|\Omega\|(C+C^2+\cdots+C^{d-2n})$ and  put $\eta=\~\eta+\eta'$ so
that $\Omega=\xi\land\eta+\Theta$. Since
$\|\~\eta\|\le\|\Omega\|$, we have
$\|\eta\|+\|\Theta\|\le\|\eta\|+\|\eta'\|+\|\Theta\|\le
(1+C+\cdots+C^{d-2n})\|\Omega\|$.
\end{proof}

\begin{Cor}\label{special-division}
If $H$ is a balanced Hamiltonian of degree $n+1$, then any
polynomial $2$-form $\Omega$ of degree $\le 3n$ can be divided by
$dH$,
\begin{equation}\label{div-corollary}
  \Omega=dH\land \eta+\Theta,\qquad
  \|\eta\|+\|\Theta\|\le(n+1)^{n+1}\cdot\|\Omega\|.
\end{equation}
\end{Cor}

\begin{proof}[Proof of the Corollary]
It is sufficient to remark that for a balanced Hamiltonian the
form $dH$ is normalized at infinity and the difference between
$dH$ and its principal homogeneous part $d\H$ is of norm $\le n$.
\end{proof}

\subsection{Derivation of the redundant Picard--Fuchs system}
Now we can write explicitly a system of first order linear
differential equations for Abelian integrals, with coefficients
explicitly bounded provided the Hamiltonian is balanced (i.e., its
lower order terms do not dominate the principal homogeneous part).
The reason why this system is called redundant, will be explained
below.

Consider $\nu=n(2n-1)$ monomial 2-forms $\Omega_i$ spanning
$\Lambda^2_{2n}$, and let $\omega_i\in\Lambda^1_{2n}$ be their
monomial primitives (arbitrary chosen), $d\omega_i=\Omega_i$, with
unit coefficients so that $\|\omega_i\|=1$ and $\|d\omega_i\|\le
2n$. Then any 2-form of degree $\le 2n$ can be represented as a
linear combination of $d\omega_i$, hence any 1-form of degree $\le
2n$ admits representation as a linear combination of $\omega_i$,
$i=1,\dots,\nu$ modulo an exact differential.

\begin{Thm}\label{main}
Let $H$ be a balanced Hamiltonian of degree $n+1$.

Then the column vector $I=(I_1(t),\dots,I_\nu(t))$ of integrals of
all monomial 1-forms $\omega_i$ of degree $\le 2n$ over any cycle
on the level curves $\{H(x,y)=t\}$ satisfies the system of linear
ordinary differential equations
\begin{equation}\label{pf}
  (t-A)\dot I=BI,\qquad I=I(t)\in\C^\nu,\
  A,B\in\Mat_{\nu\times\nu}(\C).
\end{equation}

The norms of the constant matrices $A,B$ are explicitly bounded:
\begin{equation}\label{AB-bounds}
\|A\|+\|B\|\le 6n(n+1)^{n+1}.
\end{equation}
\end{Thm}

\begin{Rem}
We use here the norms of matrices \eqref{matrix-norm}, associated
with $\ell^1$-norms on the spaces of polynomials, as defined in
\eqref{matrix-norm}.
\end{Rem}

\begin{Rem}
As was already mentioned, the assumption that $H$ is balanced,
does not involve loss of generality, since any Hamiltonian regular
at infinity can be balanced by appropriate affine transformation
(see however the discussion below).
\end{Rem}

\begin{proof}[Proof of the Theorem]
We start with a computation showing that the system can be indeed
written in the form \eqref{pf}: this derivation will be later
slightly modified to produce explicit bounds.

For any $i=1,\dots,\nu$ the 2-form $H\,d\omega_i$ of degree $\le
n+1+2n$ can be divided out with remainder by the form $dH$ (which
is, by assumption, normalized at infinity):
\begin{equation}\label{div-by-dH}
 H\,d\omega_i=dH\land\eta_i+\Theta_i,
 \qquad
 \deg\eta_i\le\deg d\omega_i\le 2n,\ \deg\Theta_i\le 2n.
\end{equation}
Since  2-forms $d\omega_i$ span the whole space of 2-forms of
degree $\le 2n$, every $d\eta_i$ and $\Theta_i$ are  linear
combinations of $d\omega_j$:
$$
 d\eta_i=\sum_{j=1}^\nu b_{ij}\,d\omega_j,
 \qquad
 \Theta_i=\sum_{j=1}^\nu a_{ij}\,d\omega_j,
$$
with appropriate complex coefficients $a_{ij},b_{ij}$ forming two
$\nu\times\nu$-matrices $A,B$ respectively, and certain
polynomials $F_j\in\C[x,y]$.

The first identity implies that $\eta_i=\sum_j
b_{ij}\omega_j+dF_i$ for suitable polynomials $F_i$. Integrating
over cycles on the level curves $\{H=t\}$ and using the
Gelfand--Leray formula for derivatives, we conclude that
$$
 t\dot I_i=\sum_j b_{ij} I_j+\sum_j a_{ij}\dot I_j,
 \qquad
 i,j=1,\dots,\nu,
$$
which is equivalent to the matrix form \eqref{pf} claimed above.

In order to place the upper bounds on the matrix norms $\|A\|$ and
$\|B\|$, we can use the bounded division lemma, but additional
efforts are required. Indeed, the normalization at infinity {\em
does not imply\/} any upper bound on the norm of the principal
part $\|\^H\|$, so the norm of the left hand side in
\eqref{div-by-dH} is apriori unbounded and does not allow for
application of \lemref{bd}.

To construct a system satisfying the inequalities
\eqref{AB-bounds}, we decompose $H$ into the principal part $\^H$
and the collection of lower terms $h=H-\H$ and treat two parts,
$\H\,d\omega$ and $h\,d\omega_i$ separately.

Let $\rho\in\Lambda^1_2$ be the 1-form $x\,dy-y\,dx$ with
$\|\rho\|=2$. Then by the Euler identity,
\begin{equation}\label{euler}
  (n+1)\H\,dx\land dy=d\H\land \rho,
\end{equation}
and therefore
\begin{equation}\label{div-principal}
  \H\,d\omega_i=\frac{d\omega_i}{dx\land dy}\cdot\H\,dx\land dy=
  \frac{d\omega_i}{(n+1)\,dx\land dy}\cdot d\H\land
  \rho=(dH-dh)\land\eta'_i,
\end{equation}
where $\|\eta'_i\|\le \|d\omega_i\|\|\rho\|/(n+1)\le 2\cdot
2n/(n+1)\le 4$. Now the term $H\,d\omega_i$ can be explicitly
expanded as
$$
 H\,d\omega_i=\^H\,d\omega_i+h\,d\omega_i
 =dH\land\eta_i'-dh\land\eta_i'+h\,d\omega_i=dH\land\eta_i'+\Omega_i',
$$
where $\|\Omega_i'\|\le \|\eta'\|\,\|dh\|+\|h\|\|d\omega_i\|\le
4n+2n= 6n$ and $\deg\Omega_i'\le 3n$. Applying
Corollary~\ref{special-division}, we write
$$
\Omega_i'=dH\land \eta_i''+\Theta_i
$$
with $\|\eta_i''\|+\|\Theta_i\|\le 6n(n+1)^{n+1}$ which together
with the previous bounds for $\|\eta_i'\|$ would imply the
inequality
\begin{equation}\label{div-by-dH-bounds}
  \|\eta_i\|+\|\Theta_i\|\le 6(n+1)^{n+2}.
\end{equation}
for the identities \eqref{div-by-dH}

Since all forms $\omega_i,d\omega_i$ are monomial with norms $\ge
1$, expanding $\eta_i$ and $\Theta_i$ leads to coefficients
satisfying the conditions
$$
 \sum_{j=1}^\mu |a_{ij}|\le\|\theta_i\|,
 \qquad
 \sum_{j=1}^\mu |b_{ij}|\le\|\Theta_i\|,
$$
which gives the required bounds on $\|A\|$ and $\|B\|$.
\end{proof}

In order to incorporate the case of quasimonic but not balanced
Hamiltonians, we derive an obvious corollary.

\begin{Cor}
If $H$ is quasimonic and the difference between $H$ and its principal
homogeneous part $\^H$ is explicitly bounded,
\begin{equation}\label{non-balance}
  \|H-\^H\|\le c,
\end{equation}
then one can choose the monomial forms so that the system
\eqref{pf} for their integrals involves the matrices $A,B$
satisfying the inequality
\begin{equation}\label{disb-bound}
  \|A\|+\|B\|\le 6(n+1)^{n+2}\cdot c^{n+1}.
\end{equation}
\end{Cor}

\begin{proof}
It is sufficient to make a transformation replacing the initial
Hamiltonian $H(x,y)$ by $c^{-(n+1)}H(cx,cy)$. This will make $H$
balanced and the main theorem applicable. Notice that such
transformation implies the change of time (the independent
variable) $t\mapsto c^{-(n+1)}$ for the resulting system
\eqref{pf}. With respect to the original variable the system
\eqref{pf} will take the form with the same matrix $B$, and $A$
multiplied by $c^{n+1}$.
\end{proof}

\begin{Rem}
Note that the system in the non-balanced case is written for forms
$\omega_i$ in general not satisfying the condition
$\|\omega_i\|=1$, as was the case with balanced Hamiltonians: the
linear rescaling $(x,y)\mapsto(cx,cy)$ results in a diagonal
transformation that is in general non-scalar on the linear space
of differential forms.
\end{Rem}

\subsection{Abelian integrals of higher degrees}
The system of differential equations \eqref{pf} holds for
integrals of the basic monomial forms $\omega_i$ generating all
polynomial differential 1-forms of degree $\le 2n$. To write an
analogous system for integrals of 1-forms of higher degrees, one
can use the fact that the integrals $\oint\omega_i$ generate the
space of all Abelian integrals as a free $\C[t]$-module, provided
that $H$ is Morse and regular at infinity, see \cite{gavrilov}.
More precisely, if $\deg\omega=d$, then for any cycle $\delta(t)$
on the level curve $\{H=t\}$ one can represent
\begin{equation}\label{gavrilov}
  \oint_{\delta(t)}\omega=\sum_{i=1}^\nu p_i(t)\oint_{\delta(t)}\omega_i,
  \qquad p_i\in\C[t],
  \quad (n+1)\deg p_i+\deg\omega_i\le\deg\omega,
\end{equation}
(in fact, it is even sufficient to take any $\mu=n^2$ forms
$\omega_i$ whose differentials span $\Lambda^2/dH\land
\Lambda^1$).

Thus the linear span of all functions $t^k I_j(t)$,
$j=1,\dots,\nu$, $0\le k\le m=\lfloor d/(n+1)\rfloor$, contains
all Abelian integrals of forms of degree $\le d$. The generators
$\{t^k I_j(t)\}_{0\le k\le m}^{1\le j\le n}$ of this system
satisfy a block upper triangular system of linear first order
differential equations obtained by derivation of \eqref{pf}:
\begin{equation}\label{block-pf}
  (t-A)\frac d{dt}(t^k I)=B\,t^k I+k(t-A)\,t^{k-1}I,
  \qquad k=1,\dots, m.
\end{equation}
This system can be written in the matrix form involving two
constant $(m+1)\nu\times (m+1)\nu$-matrices exactly as \eqref{pf}
and the entries of these matrices will be explicitly bounded,
though this time the bounds and the size of the system will depend
explicitly on $d$. Nevertheless this allows to treat integrals of
forms of arbitrary fixed degree $d$ exactly as integrals of the
basic forms.

\subsection{Properties of the redundant Picard--Fuchs system}
\label{sec:fuchs-finite} Directly from the form in which the
system \eqref{pf} was obtained, it follows that it has singular
points at all critical values $t=t_j$ of the Hamiltonian; the
eigenvector corresponding to the eigenvalue $t_j$ has coordinates
$\frac{d\omega_i}{dx\land dy}(x_j,y_j)$, $i=1,\dots,\nu$. However,
in general (since $\mu<\nu$) these eigenvalues do not exhaust the
spectrum of $A$.

The other eigenvalues of $A$ actually depend on the division with
remainder, that is non-unique because the forms $\omega_i$ are
linear dependent in $\Lambda^2/dH\land \Lambda^1$. Thus no
invariant meaning can be associated with this part of the
spectrum, and the corresponding singularities are apparent for the
Abelian integrals (though other solutions can well have
singularities at these ``redundant'' points).

However, this freedom can be used to guarantee that all these
singularities can be made Fuchsian, by slightly perturbing the
matrices.

\begin{Prop}\label{fuchs-finite}
If $H$ is a Morse function on $\C^2$, then the system \eqref{pf}
can be constructed so that the matrix $A$ has a simple spectrum
while satisfying the same inequalities as before.
\end{Prop}

\begin{Cor}
The redundant system \eqref{pf} can be always constructed having
only Fuchsian singularities on the Riemann sphere.
\end{Cor}

\begin{proof}[Proof of the Corollary]
The Fuchsian condition at infinity is satisfied automatically, as
the matrix function $(t-A)^{-1}B$ has a simple pole at $\tau=0$ in
the chart $\tau=1/t$. The inverse $(t-A)^{-1}$ can be obtained by
dividing the {\em adjugate matrix\/} (a matrix polynomial of
degree $\nu-1$ in $t$) by the determinant of $t-A$, i.e., by the
characteristic polynomial of $A$. As the latter has only simple
roots by \propref{fuchs-finite}, all poles of $(t-A)^{-1}B$ at
finite points are simple.
\end{proof}

\begin{proof}[Proof of the Proposition]
Assume that the enumeration of the forms $\omega_i$ is arranged so
that the first $\mu$ of them constitute a basis in
$\Lambda^2/dH\land \Lambda^1$. The procedure of division of the
forms $H\,d\omega_i$ by $dH$ can be altered to produce a unique
answer, if we require that the remainder is always a linear
combination of only the first $\mu$ forms. Moreover, instead of
dividing the forms $H\,d\omega_i$ with $\mu+1\le i\le\nu$, we will
divide the forms $(H-\lambda_i)\,d\omega_i$ with arbitrarily
chosen constants $\lambda_i\in\C$, $i=\mu+1,\dots,\nu$:
\begin{align*}
 H\,d\omega_i&-\sum_{j=1}^\mu a_{ij}d\omega_j\in dH\land \Lambda^1,
 \quad i=1,\dots,\mu,
 \\
 (H-\lambda_i)d\omega_i&-\sum_{j=1}^\mu a_{ij}d\omega_j\in dH\land \Lambda^1,
 \quad i=\mu+1,\dots,\nu.
\end{align*}
After division organized in such a way, the matrix $A$ of the
system \eqref{pf} obtained after expanding the incomplete
fractions, will have block lower-triangular form. The upper-left
block of size $\mu\times\mu$ has as before the eigenvalues
$t_1,\dots,t_\mu$, while the lower-right block of size
$(\nu-\mu)\times(\nu-\mu)$ is diagonal with $\lambda_i$ being the
diagonal entries. Note that in this alternative derivation we lost
control over the magnitude of the coefficients of remainders and
incomplete ratios.

Thus for the same column vector of Abelian integrals we have
constructed {\em two\/} essentially different systems of the same
form \eqref{pf} but with different pairs $(A,B)$ of
$\nu\times\nu$-matrices (the first bounded in the norm, the second
with a predefined spectrum). By linearity, any linear homotopy
between the two systems will also admit all Abelian integrals as
solutions.

Consider such a homotopy parameterized by $s\in[0,1]$. The
eigenvalues of the matrix $A$ do depend algebraically on the
parameter $s$. For $s=1$ they are equal to the critical values
$t_1,\dots,t_\mu$ of $H$ and arbitrarily prescribed values
$\lambda_{\mu+1},\dots,\lambda_\nu$. Since $t_i\ne t_j$ and
$\lambda_i$ can be also chosen different from all $t_j$ and from
each other, the eigenvalues are simple for $s=1$ and hence they
remain pairwise different for almost all values of $s$, in
particular, for arbitrarily small positive $s$ when the system is
arbitrarily close to the first system (of explicitly bounded
norm). Perturbing in that way achieves simplicity of the spectrum
of the matrix $A$ while changing the norms of $A,B$ arbitrarily
small.
\end{proof}

\section{Zeros of Abelian integrals away from the singular locus
and related problems on critical values of polynomials}
\label{motivation}

\subsection{Heuristic considerations}
As was already noted in the introduction, the main reason why so
much emphasis was put on explicit upper bounds for the
coefficients of the system \eqref{pf}, was applications to the
tangential Hilbert 16th problem. We explain in this section how
the information accumulated so far can be used to place an
effective upper bound for the number of zeros of Abelian integrals
on a positive distance from the critical (ramification) locus and
also on the total number of zeros of any branch for Hamiltonians
whose critical values are distant from each other.

However, the group of affine transformations acts naturally on the
space of Abelian integrals in a quasihomogeneous manner, and to be
geometrically sound, upper bounds for the number of zeros should
be compatible with this symmetry. In particular, the above
mentioned ``positive distance to the critical locus'' (resp.,
``distance between the critical values'') should be invariant by
affine rescaling of Hamiltonians. Besides intrinsic
considerations, the need for the bounds invariant by this action
is motivated by the future study of zeros of Abelian integrals
near singularities (cf.~with \cite{era-99}).

>From the analytic point of view, the problem is in the choice of
normalization on the variety of Hamiltonians regular at infinity.
The geometric invariance requires this normalization to be imposed
in terms of geometry of {\em configurations of the critical
values\/} of the Hamiltonians. On the other hand, the assertion of
the theorem on zeros of Abelian integrals, derived from the
explicit form of the system \eqref{pf}, uses pre-normalization in
terms of the {\em coefficients\/} of the Hamiltonian, more
precisely, the $\ell^1$-norms of its nonhomogeneity (the
difference between $H$ and its principal homogeneous part).

Thus in a natural way the problem on equivalence of the two
normalizing conditions arises. It can be shown relatively easily
that a quasimonic polynomial whose non-principal part is bounded
from above (in the sense of the norm), has all critical values
inside a disk of known radius shrinking to a point as the
non-principal part tends to zero (\propref{inverse-sense} below).

One might hope that a converse statement is also true: if all
critical values of a Hamiltonian $H$ come very close to each
other, then (eventually after appropriate translations in the
preimage and the image) $H$ differs from its principal homogeneous
part $\H$ by a small polynomial.

This fact indeed holds true for univariate (and hence
hyperelliptic) polynomials, where we were able to produce explicit
inequalities between the diameter of the critical locus
$\operatorname{diam}\S=\max_{i,j=1,\dots,\mu}|t_i-t_j|$ and the
nonhomogeneity  $\|H-\H\|$, see \thmref{roots} and
Corollary~\ref{cor2roots}.

Yet for the truly bivariate polynomials the problem turned out to
be considerably harder, and the best we were able to do is to show
that {\em for any fixed principal part $\H$ the above two
normalizations are equivalent\/}, but as different linear factors
of $\H$ approach each other, the equivalence explodes.

The current section explains the above arguments in more details
and introduces the problem on relationships between the spread of
critical values and an effective nonhomogeneoty of polynomials.
Known partial results in this sense are collected in the next
section.

\subsection{Meandering theorem and upper bounds for zeros of
Abelian integrals} A (scalar) linear ordinary differential
equation with explicitly bounded coefficients admits an explicit
upper bound for the number of isolated (real or complex) zeros of
all its solutions, see \cite{jde-95,fields}.

The {\em system\/} of equations \eqref{pf} can be reduced to one
linear equation of degree $\le \nu^2$ with rational in $t$
coefficients in such a way that any linear combination
$u(t)=\sum_{i=1}^\nu c_i I_i(t)$ of the integrals $I_i(t)$  with
constant coefficients $c_1,\dots,c_\nu\in\C$, will be a solution
to this equation: it is sufficient to find a linear dependence
between any fundamental $\nu\times\nu$-matrix $X(t)$ and its
derivatives up to order $\nu^2-1$ over the field $\C(t)$ of
rational functions.

Unfortunately, this procedure does not allow to place any bound on
the magnitude of coefficients of the resulting equation. Instead,
in \cite[Appendix B]{meandering} we described an algorithm of
derivation of another linear equation of much higher order, whose
coefficients are polynomially depending on the coefficients of the
initial system \eqref{pf}. This algorithm is explicit, so that all
degrees and coefficients admit explicit upper bounds. As a result,
the system \eqref{pf} is reduced to a {\em Fuchsian linear
differential equation\/} of the form
\begin{equation}\label{pf-eqn}
\begin{gathered}
 \Delta^\ell(t)\,u^{(\ell)}+h_{\ell-1}(t)\Delta^{\ell-1}(t)\,u^{(\ell-1)}+
 \cdots+h_1(t)\Delta(t)\,u'+h_0(t)\,u=0,
  \end{gathered}
\end{equation}
where $\Delta(t)=(t-t_1)\dots(t-t_\nu)$ is the characteristic
polynomial of the matrix $A$ and all polynomial coefficients
$h_i(t)\in\C[t]$, $i=0,\dots,\ell-1$, have degrees $\deg h_i$ and
heights $\|h_i\|$ explicitly bounded by elementary functions of
$n$. It is important to note here that the bounds, though
completely explicit, are enormously excessive, being towers
(iterated exponents) of height $4$.

The coefficients of the equation \eqref{pf-eqn} are explicitly
bounded from above on the complement to sublevel sets
$\{|\Delta(t)|\ge \e\}$ for every given positive $\e>0$. At the
roots of $\Delta$ (eigenvalues of $A$) the equation \eqref{pf-eqn}
has Fuchsian singularities, but the eigenvalues of $A$ that are
{\em not\/} critical values of $H$, are {\em apparent
singularities\/} for all linear combinations of the Abelian
integrals $I_j$ (see \secref{fuchs-finite}).

Recall that $\S$ is the critical locus (collection of all critical
values) of the Hamiltonian $H$. Let $R$ be a finite positive
number and $K_R\Subset\C\ssm\S$ the set obtained by cutting the
set
\begin{equation}\label{KR}
 \{t\in\C\: \forall j=1,\dots,\mu\ |t-t_j|> 1/R,\ |t|< R\}
\end{equation}
along no more than $\mu$ line segments to produce a {\em simply
connected\/} compact ``on the distance $1/R$ from both $\S$ and
infinity''.

Applying a general theorem on oscillations of solutions of linear
equations with bounded coefficients
\cite{lleida,meandering,fields}, we arrive to the following
theorem.

\begin{Thm}[see \cite{meandering}]\label{onzeros}
Let $H$ be a balanced Hamiltonian of degree $n+1$ and $K_R$ a
compact on distance $1/R$ from the critical locus of $H$ in the
sense of \eqref{KR}.

Then the number of zeros inside $K_R$ of any Abelian integral of a
form of degree $d$ does not exceed $(2+R)^N$, where $N=N(n,d)$ is
a certain elementary function depending only on $n$ and $d$.
\end{Thm}

The function $N(n,d)$ can be estimated from above by a tower of
four stories (iterated exponent) and certainly gives a very
excessive bound. Yet we would like to remark that this is
absolutely explicit bound, involving no undefined constants.

\begin{Rem}
The necessity of cutting in the definition of $K_R$ is due to the
fact that Abelian integrals are multivalued and a choice of branch
should be specified each time when zeros are counted.
\end{Rem}

The coefficients of the equation \eqref{pf-eqn} blow up as
$t\to\Sigma$, so no upper bound for zeros can be derived from the
general theorem \cite{jde-95}. However, if the singularity $t_i$
is apparent and distant from all other points, say, at least by
$1$, then one can place an upper bound on the coefficients of
\eqref{pf-eqn} on the boundary of the disk $\{|t-t_i|=\tfrac12\}$
and then by \cite[Corollary 2.7]{fields} the variation of argument
of any solution along the boundary can be explicitly bounded and
by the argument principle, this would imply an upper bound for the
number of zeros also {\em inside\/} the disk, where the
coefficients are very large.

It turns out that a similar construction can be also carried out when
$t_i$ is a true (non-apparent) singularity, provided that  it is of
Fuchsian type and the spectrum of the monodromy operator is on the unit
circle.

Suppose that a function $u(t)$ analytic in the punctured disk
$\{0<|t-t_i|\le 1\}$ admits a finite representation
$u(t)=\sum_{\lambda,k}f_{k,\lambda}(t)(t-t_i)^\lambda
\ln^k(t-t_i)$ with coefficients $f_{k,\lambda}$ analytic in the
closed disk $\{|t-t_i|\le 1\}$, involving {\em only real\/}
exponents $\lambda$. If this function satisfies a linear ordinary
differential equation (with a Fuchsian singularity at $t=t_i$)
whose coefficients are explicitly bounded on the boundary
circumference of this disk, then it is proved in \cite[Theorem
4.1]{fields} that any branch of $u$ admits an upper bound for the
number of zeros in this disk in terms of the magnitude of the
coefficients on the boundary and the order of the equation (the
first result of this type was proved in \cite{roitman}).

The assumption on the spectrum always holds for Abelian integrals,
since the above exponents $\lambda$ are always rational \cite{avg}
(in particular, equal to $1$ for a Morse critical value). Thus the
above result (together with the bounded meandering principle) can
be applied to the tangential Hilbert problem provided that all
critical values of the Hamiltonian are at least $1$-distant from
each other. An arbitrary Morse Hamiltonian one can rescaled to
such form, yet the number $\min_{i\ne j}|t_i-t_j|$ will enter then
into the expression for the bound.

By analogy with the previous result, denote by $K_\infty$ a simply
connected open set obtained by slitting $\C\ssm\S$ along rays
connecting critical values with infinity.

\begin{Thm}\label{noncollision}
Let $H$ be a balanced Hamiltonian of degree $n+1$, whose critical
values $t_1,\dots,t_\mu$ satisfy for some positive $R<\infty$ the
condition
$$
 |t_i-t_j|\ge 1/R,\quad |t_i|\le R\qquad \forall i\ne j.
$$

Then the number of zeros inside $K_\infty$ of any Abelian integral
of a form of degree $d$ does not exceed $(2+R)^{N'}$, where
$N'=N'(n,d)$ is a certain elementary function depending only on
$n$ and $d$.
\end{Thm}

\begin{proof}[Sketch of the proof]
Multiplying the Hamiltonian by $R$ and applying the bounded
meandering principle to the Picard--Fuchs system \eqref{pf}, we
construct a scalar linear equation of a very large order,
satisfied by all Abelian integrals, so that its coefficients are
explicitly bounded on distance $\ge 1$ from the critical locus by
an expression polynomial in $R$ as above.

To count zeros of Abelian integral inside the set $K_{\frac12}$,
one can use \thmref{onzeros}. The remaining part $K_\infty\ssm
K_{\frac12}$ consists of {\em disjoint\/} disks of radius $1/2$
centered at the critical values $t_i$ and slit along radii.
Theorem~4.1 from \cite{fields} applies to ever such disk and gives
an upper bound for the number of zeros in these disks, thus
completing the proof.
\end{proof}

Note the difference between two apparently similar results:
\thmref{onzeros} gives a {\em uniform\/} upper bound for the
number of zeros in a {\em certain domain\/} (depending on the
Hamiltonian, but always nonvoid for Morse Hamiltonians regular at
infinity).

On the contrary, \thmref{noncollision} formally solves the
tangential Hilbert problem for {\em all\/} Morse Hamiltonians
(giving an upper bound for the number of {\em all\/} zeros,
wherever they occur), but the bound is {\em not uniform\/} and
explodes when the Hamiltonian approaches the boundary of the set
of Morse polynomials regular at infinity.

\subsection{Affine group action and equivariant problem on zeros of
Abelian integrals} Consider the affine complex space of
Hamiltonians $\mathcal H=\Lambda^0_{n+1}$ and the space of 1-forms
$\mathcal F=\Lambda^1_d$ of a given degree $d$. The Abelian
integrals are multivalued functions on $\bigl((\C\times\mathcal
H)\ssm\boldsymbol\S\bigr)\times\mathcal F$, where $\boldsymbol\S$
is the {\em global discriminant\/},
$$
 \boldsymbol\S\subset\C\times\mathcal H,
 \qquad
 \boldsymbol\S=\{(t,H)\:\text{ $t$ is a critical value of $H$}\}.
$$
The group $G_2$ of affine transformations of $\C^2$ and the group
$G_1$ of affine transformations of $\C^1$ act naturally on
$\C\times\mathcal H$,
$$
 (H,t)\overset{g_2,g_1}{\longmapsto}
 (g_1\circ H\circ g_2, g_1 t),
$$
leaving $\boldsymbol\S$ invariant. The problem of counting zeros
of Abelian integrals should be also formulated for subsets in
$(\C\times\mathcal H)\ssm\boldsymbol\S$ that are invariant by this
action.

To achieve this equivariant formulation, we follow the ideology of
normal forms and choose a convenient representative from each
orbit of the group action. To factorize by the action of $G_1$, we
notice that any point set $t_1,\dots,t_\mu$ not reducible to one
point, can be put by a suitable affine transformation (or, what is
equivalent, by the choice of a chart) to a configuration
satisfying two conditions,
\begin{equation}\label{norm-sigma}
  t_1+\cdots+t_\mu=0,\qquad \max_{t=1,\dots,\mu}|t_j|=1,
\end{equation}
and such transformation is determined uniquely modulo rotation of
$\C$, preserving the Euclidean metric on $\C\simeq\R^2$. Any set
$\mathcal K$ in $(\C\times \mathcal H)\ssm\boldsymbol\S$ invariant
by the $G_1$-action, leaves its trace on the $t$-plane as a subset
$K$ disjoint from the points $\S=\{t_j\}_1^\mu$ and the distance
from $K$ to $\S$ measured in this privileged chart, is the natural
equivariant distance between $\mathcal K$ and $\boldsymbol\S$.

We arrive thus to the following equivariant formulation of the
problem on zeros of Abelian integrals, restricted in the sense
that it concerns only zeros distant from singularities (this
terminology was recently suggested by Yu.~Ilyashenko).

\begin{Prob}[Equivariant restricted tangential Hilbert 16th
problem]\label{equivariant} Let $H$ be a Hamiltonian of degree
$n+1$ regular at infinity, whose critical values
$t_1,\dots,t_\mu$, $\mu=n^2$, satisfy the normalizing conditions
\eqref{norm-sigma}.

For any finite  $R>0$ it is required to place an upper bound for
the number of isolated zeros of Abelian integrals
$\oint_{H=t}\omega$ of any form of degree $\le d$ in  the sets
$K_R$ as in \eqref{KR}. The bound should depend only on $n,d$ and
$R$.
\end{Prob}

\subsection{From \thmref{onzeros} to Equivariant problem}
In order to derive from \thmref{onzeros} a solution to the
equivariant problem, one should try to find in the orbit of the
$G_2$-action on $\mathcal H$ a Hamiltonian as close to be balanced
as possible.

Indeed, if for some affine transformation $g\in G_2$ the
Hamiltonian $\~H=H\circ g$ is already balanced, then integrals of
any form $\omega$ over any level curve $H=t$ are equal to
integrals of the form $g^*\omega$ over the curve $\~H=t$ (by the
simple change of variables in the integral). But as
$g\:\C^2\to\C^2$ is an affine map, the form $g^*\omega$ is again a
polynomial 1-form of the same degree as $\omega$, while the new
Hamiltonian $\~H$ is balanced. Hence \thmref{onzeros} can be
applied to produce the upper bound for the number of zeros exactly
in the form we need to solve the equivariant problem: the result
will be automatically a bound polynomial in $R$ with the exponent
depending only on $d$ and $n$.

In fact, it is sufficient to find in the $G_2$-orbit of $H$ a
Hamiltonian $\~H$ that would be quasimonic and whose difference
from its principal homogeneous part $\H$ would be of norm
explicitly bounded in terms of $n$. Indeed, if $\~H$ is such a
polynomial and $\|\~H-\H\|\le \tau=\tau(n)$, then the
transformation
\begin{equation}\label{restorebalance}
  \~H(x,y)\rightsquigarrow H^*(x,y)=\tau^{-(n+1)}\~H(\tau x,\tau y)
\end{equation}
will preserve the principal homogeneous part $\H$ while dividing
all other terms by appropriate positive powers of $\tau$ so that
in any case $\|H^*-\H\|\le 1$. This means that $H^*$ is balanced
and \thmref{onzeros} can be applied and will give a bound on zeros
$1/R$-distant from the critical locus of $H^*$ in terms of $R,n,d$
as required. The transformation \eqref{restorebalance} does not
preserve the normalizing conditions \eqref{norm-sigma}, but the
{\em conclusion\/} of \thmref{onzeros} can be rescaled to produce
an upper bound on zeros $1/\tau^{n+1}R$-distant from the
(normalized) critical locus of $H$, by a suitable power of
$R\tau^{n+1}$, which will give a solution to the equivariant
problem.

Recall that the balance condition consists of the two parts: the
(quasimonic) normalization of the principal homogeneous terms and
the unit bound for the norm of all non-principal terms. The first
part can be easily achieved by a suitable $G_2$-action. Indeed,
replacing $H(x,y)$ by $H(\tau x,\tau y)$, one can effectively
multiply the principal homogeneous part $\H$ by $\tau^{n+1}$ and
thus achieve the required normalization.

It will be convenient in the future not to change the principal
part any more, once it was made quasimonic. This means that the
only remaining degree of freedom to use is the group of
translations of $\C^2$ (and rotations that do not affect norms).

Summarizing this discussion, we see that in order to derive from
\thmref{onzeros} the equivariant restricted tangential Hilbert
16th problem (\probref{equivariant}), it would be sufficient to
solve the following problem.

\begin{Def}
For a quasimonic polynomial $H$ with the principal part $\H$ we
call its {\em effective nonhomogeneity\/} the lower bound
\begin{equation}\label{eff-non}
  \varkappa(H)=\inf_{T\in G_2}\|H\circ T-\H\|,
  \qquad T\text{ a translation of }\C^2.
\end{equation}
\end{Def}

\begin{Prob}\label{bound-eff-non}
Given a quasimonic Hamiltonian $H$ of degree $n+1$, whose critical
values satisfy the normalizing conditions \eqref{norm-sigma},
place an upper bound for the effective nonhomogeneity
$\varkappa(H)$.
\end{Prob}

This and related problem, completely independent from all previous
considerations, is discussed and partially solved in the next
section.

\section{Critical values of polynomials}\label{spread}

\subsection{Geometric consequences of quasimonicity}
The normalizing condition at infinity (for 1-forms and
Hamiltonians) was introduced in purely algebraic terms as an
inequality imposed on the principal homogeneous part of a 1-form
(resp., Hamiltonian). However, one can provide a simple geometric
meaning to this condition.

Recall that if $H$ is regular at infinity, then its principal
homogeneous part $\H$ has an isolated critical point at the
origin. This means that the gradient  $\nabla \H$ never vanishes
outside the origin, and in particular its minimal (Hermitian)
length on the boundary of the unit bidisk $\mathbb B=\{|x|\le 1,\
|y|\le 1\}\subset\C^2$ is strictly positive. Because of the
homogeneity, this is sufficient to place a lower bound on the
length of $\nabla \H$ everywhere on $\C^2\ssm\{0\}$.

\begin{Prop}\label{geom-quasimonic}
If $\H$ is normalized \rom(quasimonic\rom), then everywhere on the
boundary of the unit bidisk $\mathbb B$ the Hermitian length of
$\nabla\H$ is no smaller than $1$.
\end{Prop}

\begin{proof}
Consider the part $\partial\mathbb B_1$ of the boundary
$\partial\mathbb B$ which is given by the inequalities $|x|=1$,
$|y|\le 1$ (the other part is treated similarly). The homogeneous
polynomial $x^{2n-1}$ can be represented as $a\H_x+b \H_y$ with
$\|a\|+\|b\|\le 1$. Restricting this on $\partial\mathbb B_1$ we
see that the Hermitian product of the gradient
$\nabla\H=(\H_x,\H_y)$ and the vector field $V$ with coordinates
$(\bar a,\bar b)$ is everywhere equal to 1 in the absolute value.
The Hermitian length of $V$ at any point of $\partial\mathbb B_1$
can be easily majorized by $\sqrt{|\bar a(x,y)|^2+|\bar
b(x,y)|^2}$ which is no greater than $\sqrt{\|a\|^2+\|b\|^2}\le 1$
on $\partial\mathbb B_1=\{|x|=1,\,|y|\le 1\}$. But then by the
Cauchy inequality, the length of $\nabla\H$ cannot be smaller than
$1$ on $\partial\mathbb B_1$.
\end{proof}

\subsection{Almost-homogeneity implies close critical values}
We begin by showing that a quasimonic Hamiltonian whose
non-principal part is bounded, admits an upper bound for the
moduli of critical values. This solves the problem inverse to
\probref{bound-eff-non}.

\begin{Prop}\label{inverse-sense}
If $H$ is a quasimonic Hamiltonian of degree $n+1$ with the principal part
$\H$ and $\|H-\H\|\le\frac1{n\sqrt2}$, then the critical values of $H$ are
all in the disk $\{|t|\le 3/n\}$.
\end{Prop}

\begin{proof}
Denote $H=\H+h$. The gradient of each monomial of degree $\le n$
has the Hermitian length bounded by $n\sqrt 2$ on the unit bidisk
$\mathbb B$. Thus if $\|h\|<\frac1{n\sqrt2}$, then $\nabla h$ has
its length strictly bounded by $1$ everywhere in $\mathbb B$.

By \propref{geom-quasimonic}, the length of $\nabla\H$ is at least
$1$ everywhere on the boundary of $\mathbb B$, so by the
topological index theorem, all $\mu=n^2$ critical points of $H$
must be be inside $\mathbb B$.

Note that a quasimonic principal part $\H$ admits no apriori upper bound
on $\mathbb B$, however, the critical values of $H=\H+h$ can be explicitly
majorized. Indeed,  at any critical point $(x_*,y_*)$, the gradient of $H$
vanishes so $\nabla\H(x_*,y_*)=-\nabla h(x_*,y_*)$. By the Euler identity,
$(n+1)\,|\H(x_*,y_*)|=|(x_*,y_*)\cdot\nabla\H(x_*,y_*)|=
|(x_*,y_*)\cdot\nabla h(x_*,y_*)|\le \sqrt 2$, since the Hermitian length
of $\nabla h(x,y)$ is explicitly bounded by $1$ in $\mathbb B$. Finally,
since $|h(x,y)|\le \|h\|\le \frac1{n\sqrt 2}$ in $\mathbb B$, we conclude
that $|H(x_*,y_*)|\le |\H(x_*,y_*)|+|h(x_*,y_*)|\le \frac{\sqrt
2}{n+1}+\frac1{n\sqrt 2}\le 3/n\sqrt 2\le 3/n$.
\end{proof}

Thus when discussing the equivariant restricted Hilbert problem,
only the other direction (\probref{bound-eff-non}) is interesting.

\subsection{Dual formulation, limit and existential problems}
\probref{bound-eff-non} can be reformulated in dual terms as
follows.

\begin{Prob}[dual to \probref{bound-eff-non}]\label{dual}
Given a quasimonic Hamiltonian $H$ of effective nonhomogeneity
$\varkappa(H)=1$, place a {\em lower\/} bound on the diameter of
its critical values
$$
 \diam\S=\max_{1\le i\ne j\le \mu}|t_i-t_j|.
$$
\end{Prob}

Having solved this problem, one can easily derive from it by the
rescaling arguments as above a solution to \probref{bound-eff-non}
and vice versa.

The dual formulation of \probref{dual} allows a limit version: one
is required to show that if $H$ cannot be reduced to a homogeneous
polynomial by a translation, i.e., $\varkappa(H)>0$, then
$\diam\S>0$, i.e., not all critical values coincide.

This limit problem can be settled.

\begin{Thm}\label{one-val}
If a polynomial  $H(x,y)$ regular at infinity has only one
critical value \rom(necessarily of multiplicity $\mu=n^2$\rom),
then by a suitable translations in the preimage and the image $H$
can be made homogeneous: $H(x,y)=\H(x+\alpha,y+\beta)+\gamma$,
where $\H$ is the principal homogeneous part of $H$.
\end{Thm}

We postpone the proof of \thmref{one-val}, deriving first as a
corollary an existential solution of either of the two equivalent
Problems \ref{bound-eff-non}~and~\ref{dual}.

\begin{Cor}\label{existential}
For a quasimonic Hamiltonian $H=\H+h$ of degree $n+1$ there exist
two positive finite constants, $\alpha=\alpha(\H)$ and
$\beta=\beta(\H)$, depending only on the principal part $\H$, such
that the critical locus $\S=\S(H)$ and the effective
non-homogeneity $\varkappa(H)$ are related as follows:
\begin{equation}\label{existential-inequalities}
 \begin{aligned}
 \varkappa(H)\ge 1&\implies\S\cap\{|t|>\alpha\}\ne\varnothing,
 \\
 \S\subseteq\{|t|\le 1\}&\implies\varkappa(H)\le\beta.
 \end{aligned}
\end{equation}
\end{Cor}

\begin{proof}[Proof of the Corollary]
Consider the affine space $\mathcal H_n\simeq\C^{(n+1)(n+2)/2}$ of
polynomials of degree $\le n$ in two variables, and define two
nonnegative functions on it,
$$
 f(h)=\varkappa(\H+h)=\inf_{T\in\C^2}\|T^*(\H+h)\|,
 \qquad g(h)=\sum_{t\in\S(\H+h)}|t|^2,
$$
where $T$ ranges over all translations of the plane $T^2$ and
$\S(H)$ is the collection of all critical values of the
Hamiltonian $H=\H+h$ with the fixed principal part $\H$.

Both functions, as one can easily see, are semilagebraic on
$\C^2\simeq\R^4$. From \thmref{one-val} it follows that $f(h)$
must vanish if $g(h)=0$, i.e., that the zero locus of $g$ is
contained in that of $f$.

By the \polishL ojasiewicz inequality, there exist two positive
finite constants $C,\rho>0$, such that
$$
f(h)\le Cg^\rho(h),\qquad\forall h\in\mathcal H_n.
$$
From this inequality the assertion of the Corollary easily follows
if we let $\beta=Cn^\rho$ and $\alpha=(nC)^{-1/2\rho}$. Since
$C,\rho$ depend only on the construction of $f,g$, that is, on
$\H$, the Corollary is proved.
\end{proof}

Unfortunately, the proof gives no means to compute explicitly the
bounds $\alpha$ and $\beta$. Moreover, below we will show that
they {\em cannot be chosen uniformly\/} over all quasimonic
principal parts.

\subsection{Parallel problems for univariate polynomials}
One can easily formulate analogs of all the above problems for
univariate polynomials, in which case monic rather than quasimonic
polynomials are to be considered. Note that the critical values of
the hyperelliptic Hamiltonian $H(x,y)=y^2+p(x)$ coincide with that
of the univariate potential $p\in\C[x]$, and also the effective
nonhomogeneity (more accurately, non-quasihomogeneity) of $H$
coincides with that $\varkappa(p)$. Thus all results proved below,
are valid not only for univariate polynomials, but also for
hyperelliptic bivariate Hamiltonians.

The limit problem for this case is fairly elementary. It was
solved by A.~Chademan \cite{chademan} as a step towards the
existential solution of \probref{bound-eff-non} for univariate
polynomials, see Corollary~\ref{cor-chademan} below.

\begin{Prop}[Chademan \cite{chademan}]
A complex polynomial that has only one critical value at $t=0$, is
a translated monomial $\alpha(x-a)^{n+1}$.
\end{Prop}

\begin{proof}
Assuming without loss of generality that the polynomial $p(x)$ is
monic, we can always write the derivative
$$
p'(x)=(n+1)(x-a_1)^{\nu_1}\cdots(x-a_k)^{\nu_k},
$$
where $a_1,\dots,a_k$ are geometrically distinct critical points
and all $\nu_k>0$. For any $j=1,\dots,k$ the polynomial $p$ can be
expressed as the primitive of $p'$ integrated from $a_j$,
$$
 p(x)=p(a_j)+\int_{a_j}^xp'(s)\,ds=0+(x-a_j)^{\nu_j+1}q_j(x),
 \qquad q_j\in\C[x],
$$
in other words, $p$ is divisible by $(x-a_j)^{\nu_j+1}$. As this
holds for all points $a_j$, $j=1,\dots,k$, hence $\deg p\ge\deg
p'+k$ and therefore only one $\nu_j$ can be different from zero.
\end{proof}

In the standard way (see the demonstration of
Corollary~\ref{existential} above) the following corollary can be
derived.

\begin{Cor}[A.~Chademan \cite{chademan}]\label{cor-chademan}
If $p(x)=x^{n+1}+p_{n-1}x^{n-1}+\cdots+p_1 x+p_0$ is a monic
polynomial of degree $n+1$ without the term $x^n$, and all complex
critical values of $p$ lie in the unit disk $\{|t|\le 1\}$, then
$$
 |p_{n-1}|+\cdots+|p_1|+|p_0|\le C_n,
$$
where $C_n$ is a constant depending only on $n$.\qed
\end{Cor}

However, in the same way as before, the proof based on solution of
the limit problem gives no possibility of effectively computing
the constant $C_n$. We compute it using alternative approach.

\subsection{Spread of roots {\em vs.}~spread of critical values for
univariate monic complex polynomials}

\begin{Thm}\label{roots}
If all critical values $\{t_1,\dots,t_n\}$ of a monic univariate
polynomial $p(x)=\prod_{j=0}^n(x-x_j)$ are in the unit disk, then
the diameter of the set of its roots is no greater than $4e$\rom:
$$
 \S\subset\{|t|\le 1\}
 \implies \forall j,k=0,\dots,n\ |x_j-x_k|\le 4e.
$$
\end{Thm}

\begin{proof}
Consider the real-valued function $f\:\C\to\R$, $f(x)=|p(x)|$. It
is smooth outside the roots of the polynomial $p$. Moreover, its
critical values (different from zero) coincide with $|t_j|$, as
the critical points for $f$ and $p$ are the same.

By the main principle of the Morse theory, all sublevel sets
$M_s=\{x\in\C\:f(x)\le s\}$ for $0<s<\infty$ of the function $f$
remain homeomorphic to each other until $s$ passes through a
critical value of $f$. One can easily verify that the $M_s$ is
simply connected for all large $s$ (it differs only slightly from
the disk $\{|t|\le s^{1/n+1}\}$). Our assumption on the critical
values guarantees that the set $M_1=\{|p(x)|\le 1\}\subset\C$
corresponding to $s=1$ is therefore also connected (though its
shape can be very non-circular anymore).

On the other hand, by the famous Cartan lemma \cite{levin} for any
positive $\e$ one can delete from $\C$ one or several disks with
the sum of diameters less than $\e$ so that on the complement the
monic polynomial of degree $n+1$ satisfies a lower bound
$|p(x)|\ge (\e/4e)^{n+1}$. This lemma implies that the set $M_1$
can be covered by one or several circular disks with the sum of
diameters $\le 4e$.

But the set $M_1$ (like all sets $M_s$ with positive $s$) contains
all roots of $p$, so if there are two roots $x_i,x_j$ on the
distance more than $4e$, then the union of disks covering these
two roots simultaneously, cannot be connected (it is sufficient to
project all the disks on the line connecting these roots and
reduce the assertion to one dimension). This contradiction proves
the theorem.
\end{proof}

\begin{Cor}\label{cor2roots}
By a suitable translation $p(x)\mapsto p(x+a)$ a monic polynomial
$p(x)=x^{n+1}+\cdots$ whose critical values are normalized by the
conditions \eqref{norm-sigma}, can be reduced to the form
$p(x)=x^{n+1}+\sum_{j=0}^{n} p_j x^j$ with $\sum_j
|p_j|=\|x^{n+1}-p\|\le 8^{n+1}$.
\end{Cor}

\begin{proof}
By \thmref{roots}, the roots of $p$ form a point set of diameter $d\le 4e$
in the $x$-plane. Any such set can  be covered by a regular hexagon with
the opposite sides being at the distance $d$
\cite{bonnesen-fenchel,grunbaum}. Shifting the origin at the center of
this hexagon makes all roots $x_j$ satisfying the inequality $|x_j|\le
d/\sqrt3$.

A monic polynomial of degree $n+1$ with all roots inside the disk
of radius $r>0$ has all its coefficients bounded by the respective
coefficients of the polynomial $(x+r)^{n+1}$, by the Vieta
formulas. For the latter polynomial the sum of (absolute values
of) all coefficients is the value at $x=1$ (since all these
coefficients are nonnegative). Putting everything together, we
conclude that after shifting the origin at the center of the
hexagon, $\|p(x)\|\le (1+4e/\sqrt 3)^{n+1}\le 8^{n+1}$.
\end{proof}

\begin{Rem}
Simply shifting the origin to one of the roots makes all of them
being in the circle of radius $4e$, which finally yields an upper
bound $\|p\|\le (1+4e)^{n+1}\le 12^{n+1}$ without referring to the
claim on hexagonal cover.
\end{Rem}

The assertion of \thmref{roots} for {\em real\/} polynomials
having {\em only real\/} critical points, can be proved in a
completely different way. The following proposition gives an
insight as to how accurate the bound established in \thmref{roots}
is.

\begin{Prop}\label{chebyshev}
A monic real polynomial of degree $n+1$ with all critical points
real and all critical values in the interval $[-1,1]$, has all its
real roots in some interval of the length $4$.
\end{Prop}

\begin{proof}
Between any two roots the polynomial satisfies the condition
$-1\le p(x)\le 1$, since all critical values lie on that interval.

Among monic polynomials of degree $n+1$ on the unit interval
$-1\le x\le 1$ the smallest uniform upper bound $c_n=2^{-(n+1)}$
is achieved for the Chebyshev polynomial $T_n(x)=2^{-(n+1)}\cos
(n+1)\arccos x$: for any other monic polynomial of this degree,
the $C^0$-norm $\max_{-1\le x\le 1}|p(x)|$ will be greater or
equal to $c_n$. Applying this assertion to the polynomial
$2^{n+1}p(x/2)$ we conclude that the largest real interval on
which the monic polynomial can satisfy the condition $|p|\le 1$,
is of length $4$ (twice the length of $[-1,1]$).
\end{proof}

Thus \thmref{roots} can be considered as generalizing (in some
sense) the extremal property of the Chebyshev polynomials to the
complex domain.

\subsection{Demonstration of \thmref{one-val}}\label{sec:demo1val}
The proof of \thmref{one-val} is an immediate corollary to the two
following lemmas.

\begin{Lem}\label{one-point}
A polynomial regular at infinity and having only one complex
critical value, has a unique critical point.
\end{Lem}

This lemma is in fact valid for polynomials of any number of
variables. The second claim is dimension-specific.

\begin{Lem}\label{homogen}
A bivariate polynomial regular at infinity and having a unique
complex critical point at the origin, is homogeneous.
\end{Lem}

\begin{proof}[Proof of \lemref{one-point}]
Let $H_\e$ be an analytic one-parameter perturbation of the
polynomial $H_0=H$, such that for all $\e\ne0$ the polynomial
$H_\e$ is Morse.

Consider the monodromy group of the bundle $H_\e\:\C^2\to\C^1$ for
an arbitrary small $\e$. It is known \cite{avg} that vanishing
cycles form the basis of the homology of all fibers, each being a
cyclic vector (i.e., all continuations of any vanishing cycle span
the entire first homology of the typical fiber $\{H_\e=t\}$.

Suppose that there are at least two critical points $a_1\ne
a_2\in\C^2$ for $H_0$. Then for all sufficiently small $\e$ the
polynomial $H_\e$ will have two disjoint groups of critical points
with close critical values. Moreover, these groups of critical
points are well apart (say, the distance between them is never
smaller than half the distance between $a_1$ and $a_2$).

But then the vanishing cycles ``growing'' from critical points not
belonging to the same group, are also disjoint, therefore their
intersection index must be zero.

But then the Picard--Lefschetz formulas imply that the subspaces
generated by each group of vanishing cycles, must be both
invariant, which contradicts the fact that each group must consist
of cyclic elements for the monodromy.
\end{proof}

\begin{proof}[Proof of \lemref{homogen}]
Consider the one-parameter analytic (polynomial) homotopy between
$H$ and its principal part,
$H_\e(x,y)=\e^{n+1}H(\e^{-1}x,\e^{-1}y)$. Then for $\e=0$ $H_0$
coincides with the principal homogeneous part, while $H_1=H$.

The germ of $H_\e$ at the origin $x=y=0$ has a multiplicity
$\mu_\e$ (the Milnor number) that is equal to $n^2$ for any $\e$.
Indeed, by the B\'ezout theorem, the total number of critical
points of $H$ counted with multiplicities in the projective plane
$\C P^2$, is $n^2$; the condition of nondegeneracy at infinity
implies that all of them are in the finite (affine) part $\C^2$.
The uniqueness assumption means that all these $n^2$ points
coincide at the origin.

By the famous theorem due to D.~T.~L\hataccent{e} and
C.~P.~Ramanujam \cite{le-raman}, the topological type of an
analytic germ is constant along the stratum $\mu=\const$,
therefore the germs of $H_0$ and $H_1$ at the origin are
topologically equivalent, in particular, the germs of analytic
curves $\{H_0=0\}$ and $\{H_1=0\}$ in $(\C^2,0)$ are homeomorphic.

But by  the Zariski theorem \cite{zariski}, the order of a planar
analytic curve (i.e., the order of the lowest order terms which
occur in the Taylor expansion of the local equation defining this
curve) is a topological invariant. For the curve $H_0=0$ this
order is $n+1$, as the polynomial $H_0$ is homogeneous. But this
means that the lowest order of terms that may occur in $H_1$, is
also $n+1$, that is, $H_0=H_1$ and $H$ coincides in fact with its
principal homogeneous part.
\end{proof}

\begin{Rem}
Consider the gradient vector field $\nabla H$. Its principal
homogeneous part, $\nabla \H$, is a homogeneous vector field on
the plane that has an isolated singularity of multiplicity $n^2$
at the origin.

Assertion of \lemref{homogen} means that adding any nontrivial
lower order terms to $H$ would necessarily create singular points
of the gradient vector field outside the origin, thus changing the
multiplicity of what remains at the origin.

However, this assertion about {\em arbitrary\/} (not necessarily
gradient) polynomial vector fields is {\em false\/}, as the
following example shows.

\begin{Ex}[Lucy Moser--Jauslin]\label{lucy}
The nonhomogeneous vector field
$$
(x^3-y^3+x)\tfrac{\partial}{\partial
x}+(2x^3-y^3+x)\tfrac{\partial}{\partial y}
$$
has a unique singular point of the maximal multiplicity $9$ at the
origin, and the principal homogeneous part has an isolated
singularity.
\end{Ex}
\end{Rem}

\subsection{Existential bounds cannot be uniform}
As was already noted, the proof of Corollary~\ref{existential}
gives no indication on how to compute the bounds $\alpha(\H)$ and
$\beta(\H)$ for a given homogeneous part $\H$. However, the
folowing example shows that there cannot be the bound uniform over
all principal parts: as some of the linear factors approach each
other, the values of $\beta$ and $\alpha^{-1}$ may grow to
infinity.

\begin{Ex}
The form $\^H_a(x,y)=a \frac{x^{n+1}}{n+1}+\frac{y^{n+1}}{n+1}$ is
normalized for $a\ge 1$, as one can easily see by comparing the
operator of division by $d\^H=\<ax^n,y^n\>$ on $2n$-forms with
that by the ideal $\<x^n,y^n\>$.

The polynomial $H_a(x,y)=\^H_a(x,y)-x$ has critical points at
$y=0$, $x=1/\sqrt[n]{a}$ (of multiplicity $n$ for every choice of
branch of the root). The corresponding critical values all
converge to zero asymptotically as $a^{-1/n}$ as $a\to\infty$.

On the other hand, the effective nonhomogeneity of the univariate
polynomial $p_a(x)=a\frac{x^{n+1}}{n+1}-x$ (and hence the value
$\varkappa(H_a)$) remain bounded away from zero as $a\to\infty$.
Indeed, if after shifting the polynomial $p_a$ by $r=r(a)\in\C$
the coefficient before $x^n$ goes to zero, then necessarily $a
r(a)\to0$. On the other hand, the coefficient before the linear
term is equal to $1+a r(a)^n$ and hence is bounded away from zero.

Thus the bounds established in Corollary~\ref{existential}, cannot
be made uniform over all homogeneous parts. Of course, the reason
is that the space of quasimonic principal parts is not compact
(e.g., the polynomials $\H_a$ have no limit points as
$a\to\infty$). In turn, this is related to the fact that some of
the linear factors entering $\H_a$, tend to each other (as points
on the projective line $\C P^1$).
\end{Ex}

\subsection{Discussion: atypical values and singular
 perturbations}
The phenomenon occurring in the above example, might be
characteristic. When the Hamiltonian is not regular at infinity,
the Abelian integrals may have ramification points that are not
critical values of $H$. Such points, called {\em atypical
values\/}, must necessarily be singular for any system of
Picard--Fuchs equations, and are studied mostly by topological
means.

On the other hand, the fact that entries of the matrices $A,B$ may
grow to infinity as the principal part of $\H$ degenerates, means
that the system \eqref{pf} (written in the privileged chart to
make the assertion equivariant) undergoes a {\em singular
perturbation\/} (appearance of a large parameter in the right hand
side that is equivalent to putting a small parameter before some
of the higher order derivatives).

Thus we see that ``atypical singularities'' in the Picard--Fuchs
system can appear as a result of singular perturbation. The
analytic approach based on studying division by $dH$ and arguments
involving geometry of critical values, may be a complementary tool
for the study of singularities ``coming from infinity''.

\bibliographystyle{amsplain}

\end{document}